\documentclass[11pt,a4paper,twoside]{article}
\usepackage{latexsym,amsmath,amsthm,amssymb,amsfonts}
\usepackage{mathrsfs,multicol}
\usepackage{amscd}
\usepackage{cite}
\usepackage[ansinew]{inputenc}

\def\proof{{\bf Proof.}\quad}
\def\endproof{\hfill\hbox{$\sqcup$}\llap{\hbox{$\sqcap$}}\medskip}
\numberwithin{equation}{section}

\newtheorem{expl}{\indent{\bf Example}}[section]

\catcode`@=11

\newskip\plaincentering \plaincentering=0pt plus 1000pt minus 1000pt
\def\@plainlign{\tabskip=0pt\everycr={}}
\def\eqalignno#1{\displ@y \tabskip\plaincentering
  \halign to\displaywidth{\hfil$\@lign\displaystyle{##}$\tabskip\z@skip
    &$\@lign\displaystyle{{}##}$\hfil\tabskip\plaincentering
    &\llap{$\@lign##$}\tabskip\z@skip\crcr
    #1\crcr}}
\def\leqalignno#1{\displ@y \tabskip\plaincentering
  \halign to\displaywidth{\hfil$\@lign\displaystyle{##}$\tabskip\z@skip
    &$\@lign\displaystyle{{}##}$\hfil\tabskip\plaincentering
    &\kern-\displaywidth\rlap{$\@lign##$}\tabskip\displaywidth\crcr
    #1\crcr}}
\def\plainLet@{\relax\iffalse{\fi\let\\=\cr\iffalse}\fi}
\def\plainvspace@{\def\vspace##1{\noalign{\vskip##1}}}

\def\intic@{\mathchoice{\hskip5\p@}{\hskip4\p@}{\hskip4\p@}{\hskip4\p@}}
\def\negintic@
 {\mathchoice{\hskip-5\p@}{\hskip-4\p@}{\hskip-4\p@}{\hskip-4\p@}}
\def\intkern@{\mathchoice{\!\!\!}{\!\!}{\!\!}{\!\!}}
\def\intdots@{\mathchoice{\cdots}{{\cdotp}\mkern1.5mu
    {\cdotp}\mkern1.5mu{\cdotp}}{{\cdotp}\mkern1mu{\cdotp}\mkern1mu
      {\cdotp}}{{\cdotp}\mkern1mu{\cdotp}\mkern1mu{\cdotp}}}
\newcount\intno@
\def\iint{\intno@=\tw@\futurelet\next\ints@}
\def\iiint{\intno@=\thr@@\futurelet\next\ints@}
\def\iiiint{\intno@=4 \futurelet\next\ints@}
\def\idotsint{\intno@=\z@\futurelet\next\ints@}
\def\ints@{\findlimits@\ints@@}
\newif\iflimtoken@
\newif\iflimits@
\def\findlimits@{\limtoken@false\limits@false\ifx\next\limits
 \limtoken@true\limits@true\else\ifx\next\nolimits\limtoken@true\limits@false
    \fi\fi}
\def\multintlimits@{\intop\ifnum\intno@=\z@\intdots@
  \else\intkern@\fi
    \ifnum\intno@>\tw@\intop\intkern@\fi
     \ifnum\intno@>\thr@@\intop\intkern@\fi\intop}
\def\multint@{\int\ifnum\intno@=\z@\intdots@\else\intkern@\fi
   \ifnum\intno@>\tw@\int\intkern@\fi
    \ifnum\intno@>\thr@@\int\intkern@\fi\int}
\def\ints@@{\iflimtoken@\def\ints@@@{\iflimits@
   \negintic@\mathop{\intic@\multintlimits@}\limits\else
    \multint@\nolimits\fi\eat@}\else
     \def\ints@@@{\multint@\nolimits}\fi\ints@@@}
\def\Sb{_\bgroup\vspace@
        \baselineskip=\fontdimen10 \scriptfont\tw@
        \advance\baselineskip by \fontdimen12 \scriptfont\tw@
        \lineskip=\thr@@\fontdimen8 \scriptfont\thr@@
        \lineskiplimit=\thr@@\fontdimen8 \scriptfont\thr@@
        \Let@\vbox\bgroup\halign\bgroup \hfil$\scriptstyle
            {##}$\hfil\cr}
\def\endSb{\crcr\egroup\egroup\egroup}
\def\Sp{^\bgroup\vspace@
        \baselineskip=\fontdimen10 \scriptfont\tw@
        \advance\baselineskip by \fontdimen12 \scriptfont\tw@
        \lineskip=\thr@@\fontdimen8 \scriptfont\thr@@
        \lineskiplimit=\thr@@\fontdimen8 \scriptfont\thr@@
        \Let@\vbox\bgroup\halign\bgroup \hfil$\scriptstyle
            {##}$\hfil\cr}
\def\endSp{\crcr\egroup\egroup\egroup}
\def\Let@{\relax\iffalse{\fi\let\\=\cr\iffalse}\fi}
\def\vspace@{\def\vspace##1{\noalign{\vskip##1 }}}
\def\aligned{\,\vcenter\bgroup\plainvspace@\plainLet@\openup\jot\m@th\ialign
  \bgroup \strut\hfil$\displaystyle{##}$&$\displaystyle{{}##}$\hfil\crcr}
\def\endaligned{\crcr\egroup\egroup}
\def\matrix{\,\vcenter\bgroup\plainLet@\plainvspace@
    \normalbaselines
  \m@th\ialign\bgroup\hfil$##$\hfil&&\quad\hfil$##$\hfil\crcr
    \mathstrut\crcr\noalign{\kern-\baselineskip}}
\def\endmatrix{\crcr\mathstrut\crcr\noalign{\kern-\baselineskip}\egroup
                \egroup\,}
\newtoks\hashtoks@
\hashtoks@={#}
\def\format{\crcr\egroup\iffalse{\fi\ifnum`}=0 \fi\format@}
\def\format@#1\\{\def\preamble@{#1}%
  \def\c{\hfil$\the\hashtoks@$\hfil}%
  \def\r{\hfil$\the\hashtoks@$}%
  \def\l{$\the\hashtoks@$\hfil}%
  \setbox\z@=\hbox{\xdef\Preamble@{\preamble@}}\ifnum`{=0 \fi\iffalse}\fi
   \ialign\bgroup\span\Preamble@\crcr}

\def\cases{\left\{\,\vcenter\bgroup\plainvspace@
     \normalbaselines\openup\jot\m@th
      \plainLet@\ialign\bgroup$\displaystyle{##}$\hfil&
      \quad$\displaystyle{{}##}$\hfil\crcr
      \mathstrut\crcr\noalign{\kern-\baselineskip}}

\newif\iftagsleft@
\tagsleft@true
\def\TagsOnRight{\global\tagsleft@false}
\def\tag#1$${\iftagsleft@\leqno\else\eqno\fi
 \hbox{\def\pagebreak{\global\postdisplaypenalty-\@M}%
 \def\nopagebreak{\global\postdisplaypenalty\@M}\rm(#1\unskip)}%
  $$\postdisplaypenalty\z@\ignorespaces}
\interdisplaylinepenalty=\@M
\def\allowdisplaybreak{\noalign{\allowbreak}}
\def\plainallowdisplaybreak@{\def\allowdisplaybreak{\noalign{\allowbreak}}}
\def\plaindisplaybreak@{\def\displaybreak{\noalign{\break}}}
\def\align#1\endalign{\def\tag{&}\plainvspace@\plainallowdisplaybreak@
\plaindisplaybreak@
  \iftagsleft@\plainlalign@#1\endalign\else
   \plainralign@#1\endalign\fi}
\def\plainralign@#1\endalign{\displ@y\plainLet@\tabskip\plaincentering
\halign to\displaywidth
     {\hfil$\displaystyle{##}$\tabskip=\z@&$\displaystyle{{}##}$\hfil
       \tabskip=\plaincentering&\llap{\hbox{\rm(##\unskip)}}\tabskip\z@\crcr
             #1\crcr}}
\def\plainlalign@
 #1\endalign{\displ@y\plainLet@\tabskip\plaincentering\halign to \displaywidth
   {\hfil$\displaystyle{##}$\tabskip=\z@&$\displaystyle{{}##}$\hfil
   \tabskip=\plaincentering&\kern-\displaywidth
        \rlap{\hbox{\rm(##\unskip)}}\tabskip=\displaywidth\crcr
               #1\crcr}}

\def\re@#1{\par\hangindent\parindent\indent\llap{#1\enspace}\ignorespaces}
\def\qfootnote#1{\edef\@sf{\spacefactor\the\spacefactor}{}#1\@sf
      \insert\footins{\let\egroup=}\footnotesize 
      \interlinepenalty100 \let\par=\endgraf
        \leftskip=0pt \rightskip=0pt
        \splittopskip=10pt plus 1pt minus 1pt \floatingpenalty=20000
   \smallskip\re@{#1}\bgroup\strut\aftergroup{\strut\egroup}\let\next}
\catcode`\@=\active
\TagsOnRight
\begin{document}
\title{\bf The Mean Curvature Flow in Minkowski Spaces}
\author{  \ Fanqi Zeng ,\ Qun He,\ Bin Chen }
\date{}
\maketitle
\begin{quotation}
\noindent{\bf Abstract.}~ Studying the geometric flow plays a powerful role in mathematics
 and physics. In this paper, we introduce the mean curvature flow on Finsler
manifolds and give a number of examples of the mean curvature
flow. For Minkowski spaces, a special case of Finsler manifolds, we will prove the existence and uniqueness for solution of
 the mean curvature flow and prove that the flow preserves the convexity and
 mean convexity. We also derive some comparison principles for the mean curvature flow.
\end{quotation}

\section{Introduction}

Let $\varphi(\cdot, t): M^{n} \times [0, T) \to \mathbb{R}^{n+1}$ be a family of smooth closed hypersurfaces in $\mathbb{R}^{n+1}$ satisfying
\begin{equation}\label{Fin1}
\frac{\partial}{\partial t}\varphi(x, t)=\tilde{H}(x, t)\mathbf{n}(x, t),
\end{equation}
where $\tilde{H}$ is the mean curvature function and $\mathbf{n}$ is the inner pointing unit normal.
\eqref{Fin1} is the so called mean curvature flow (MCF).
Huisken \cite{Huisken84} proved that any compact convex
hypersurface in $\mathbb{R}^{n+1}$, $n \geq 2$, remains smooth and convex for a finite time under mean
curvature flow, then contracts to a point while becoming spherical in shape.

After Huisken's work in \cite{Huisken84} on mean curvature flow of convex hypersurfaces in Euclidean spaces, there have been plenty of results obtained not only for different ambient spaces, but more generally, for different
kinds of curvature flows of convex hypersurfaces in Euclidean spaces, for example, \cite{ Andrews94,Gerhardt90,Urbas90,Chow85,Schulze05 } and the references therein.
Besides  the MCF is important in many sections of mathematics and physics. It appears in applications such as models of annealing metals \cite{Mullins56} and other problems involving phase changes and moving interfaces \cite{Evans92, Ilmanen93}, and has been used to prove results in differential geometry \cite{Ilmanen86}. Mean curvature flow equations are completely difficult to be solved in all generality, due to their nonlinearity. Many aspects of mean curvature flow are well understood, including the short time existence of solutions, local regularity \cite{Ecker89}, much of the nature of singularities \cite{Huisken90}, and generalised solutions allowing continuation beyond singularities \cite{Brakke78}.

Finsler geometry is the most natural generalization of Riemannian geometry. Finsler manifolds are known to possess extremely rich geometric structures. Since the parabolic system or hyperbolic system is one of the most models in the nature, we feel the mean curvature flow on Finsler manifold is also a natural tool. As a special case of Finsler manifold, a Minkowski space is a vector space endowed with Minkowski metric which is Euclidean
metric without quadratic restriction. Curvature flows in Minkowski spaces were investigated by other authors, see e.g. \cite{Andrews, Andrews2000, Xia16}. Note that the curvature flows in \cite{Andrews, Andrews2000, Xia16} are called ``anisotropic curvature flows''. In fact, they can be thought as
curvature flows on the Minkowski spaces.

To the best of our knowledge, there are less works on higher dimensional mean curvature flow concerning
about detailed convergence on Finsler manifolds.  Compared with the case of Riemannian geometry, in our
case, we need to overcome some obstructions, the major one is how to get a priori estimate  from the PDE point of view.
Due to a lot of non-Riemannian geometric quantities, most of classical approach to prove a priori estimates by Huisken fails. Particularly, when we write the flow function as a scalar function of the graph function $r$ over the inverse Minkowski hypersphere, the evolution equation for $\partial_{t}r$ behaves not well. Also, the
evolution equation for the anisotropic second fundamental form is quite bad.

In this paper,  our motivation is to better understand the mean curvature flow on Finsler manifolds, hoping that in further studies we can detect interesting geometric objects by running the mean curvature flow or similar flows to time infinity without developing any singularity or after handling possible singularities.
As a first step, we will give a definition of the mean curvature flow on Finsler manifolds. It is not an easy way to define mean curvature flow of mutually compatible fundamental geometric structures on Finsler manifolds. Our mean curvature flow equation has the following form:
$$
\left\{\begin{array}{l}
\partial_{t}\varphi(p,t)={\hat{H}}(p,t)\textbf{n}(p,t),\\
\varphi(p,0)=\varphi_{0}(p),
\end{array}\right.
$$
where ${\hat{H}}(p,t)$ and $\textbf{n}(p,t)$ are respectively the anisotropic mean curvature and  the inner pointing unit normal
of the hypersurface $\varphi_{t}$ at the point $p\in M$.  For the exact definition of $\hat{H}$ we refer to Section 2.2.1. We also provide some examples of the mean curvature flow in Finsler geometry.
Then, due to the complexity of the evolution equation, we focus on the mean curvature flow on the Minkowski spaces in the rest of the paper and obtain some interesting results.

The contents of this paper are organized as follows.  In Section 2, we provide some fundamental concepts and formulas which are necessary for the present paper. In Section 3, we give the definition of the mean curvature flow on Finsler manifolds. In Section 4, we give some examples of the mean curvature flow in Finsler geometry. In Section 5, for Minkowski spaces we use a strictly parabolic second-order quasi-linear partial differential equation and prove the existence and uniqueness of solutions of mean curvature flow. In Section 6, we give some evolution equations. As some applications, In Section 7 and Section 8, we provide the convexity, mean convexity and comparison principles for the mean curvature flow on the Minkowski spaces.

\section{Preliminaries}

\subsection{Finsler manifolds}

Here and from now on, we will use the following convention of index
ranges unless other stated:
$$1\leq i,j,k,\cdots\leq n;\ \ \ \ \  1\leq a,b,c,\cdots\leq n-1;     \ \ \ \ \; 1\leq \alpha,\beta,\gamma,\cdots\leq n+1;\ \ \ \ \   1\leq I,J,\cdots\leq 2n$$ throughout this paper. For simplicity, from now on we will follow the summation convention and frequently
use the notations $F=F(y)$, $F_{i}=F_{y^{i}}(y)$, $u_{i}=\frac{\partial u}{\partial x^{i}}$, $u_{ij}=\frac{\partial^{2}u}{\partial x^{i}\partial x^{j}}$ and so on.

 We assume that $M$ is an $n$-dimensional
oriented smooth manifold. Let $TM$ be the tangent bundle over $M$
with local coordinates $(x,y)$, where $x=(x^1,\cdots ,x^{n})$ and
$y=(y^1,\cdots ,y^{n})$. A {\it Finsler metric} on $M$ is a function
$F: TM\longrightarrow[0,\infty)$ satisfying the following
properties: (i) $F$ is smooth on $TM\backslash\{0\}$; (ii)
$F(x,\lambda y)=\lambda F(x,y)$ for all $\lambda>0$; (iii) the
induced quadratic form $g$ is positive-definite, where
$$g:=g_{ij}(x,y)dx^{i} \otimes dx^{j}, ~~~~~~~~~   g_{ij}(x,y)=\frac{1}{2}[F^{2}] _{y^{i}y^{j}}.  $$
The projection $\pi : TM\longrightarrow M $ gives rise to the
pull-back bundle $\pi^{\ast}TM$ over $TM\backslash\{0\}$. As is well known, on  $\pi^{\ast}TM$ there
exists uniquely the \emph{Chern connection} $\nabla$ with $\nabla
\frac{\partial}{\partial x^i}=\omega_{i}^{j}\frac{\partial}{\partial
x^{j}}$ satisfying
$$\omega^{i}_{j}\wedge dx^{j}=0,$$
$$dg_{ij}-g_{ik}\omega^{k}_{j}-g_{kj}\omega^{k}_{i}=2FC_{ijk}\delta y^{k},\qquad \delta y^{k}:=\frac{1}{F}(dy^k+y^j\omega^{k}_{j}),$$
where $C_{ijk}=\frac{1}{2}\frac{\partial g_{ij}}{\partial y^k}$ is
called the\emph{ Cartan tensor}.

Let $X=X^{i}\frac{\partial}{\partial x^{i}}$ be a vector field. Then
the \emph{covariant derivative} of $X$ along
$v=v^i\frac{\partial}{\partial x^i}$ with respect to $w\in
T_{x}M\backslash \{0\}$ is defined by
$$D^{w}_{v}X(x):=\left\{v^{j}\frac{\partial X^{i}}{\partial x^{j}}(x)+\Gamma^{i}_{jk}(w)v^{j}X^{k}(x)\right
\}\frac {\partial}{\partial x^{i}},\label{Z1}$$ where
$\Gamma^{i}_{jk}=\omega^{i}_{j}(\partial x^k)$.

Let ${\mathcal L}:TM\longrightarrow T^{\ast}M$ denote the {\it
Legendre transformation}, which satisfies ${\mathcal L}(
0)=0$ and ${\mathcal L}(\lambda
y)=\lambda {\mathcal L}(y)$ for all $\lambda>0$ and $ y\in TM\setminus \{0\}$. Then ${\mathcal L}:TM\setminus
\{0\}\longrightarrow T^{\ast}M\setminus \{0\}$ is a norm-preserving
$C^\infty$ diffeomorphism.
For a smooth function $f:M\longrightarrow \mathbb{R}$, the \emph{gradient
vector} of $f$ at $x\in M$ is defined as $\nabla f(x):={\mathcal
L}^{-1}(df(x))\in T_{x}M$, which can be written as
$$\nabla f(x):=\left\{\begin{array}{l}
g^{ij}(x,\nabla f)\frac{\partial f}{\partial x^{j}}\frac{\partial }{\partial x^{i}},~~~~~~ df(x) \neq 0,\\
0,~~~~~~~~~~~~~~~~~~~~~~~~~~df(x)=0.
\end{array}\right.\label{Z2}$$
Set $M_{f}:=\{x\in M|df(x)\neq 0\}$.  We define $\nabla^{2}f(x)\in
T^{\ast}_{x}M\otimes T_{x}M$ for $x\in M_{f}$ by using the following covariant
derivative
$$\nabla^{2}f(v):=D^{\nabla f}_{v}\nabla f(x)\in T_{x}M,~~~~~~~~v\in T_{x}M.\label{Z3}$$

Let $d\mu=\sigma(x)dx^{1}\wedge\cdots\wedge dx^{n}$ be an arbitrary
volume form on $(M,F)$.  The \emph{divergence} of a smooth vector
field $V=V^{i}\frac{\partial}{\partial x^{i}}$ on $M$ with respect
to $d\mu$ is defined by
$$ \textmd{div}V :=\frac{1}{\sigma}\frac{\partial}{\partial x^{i}}\left(\sigma V^i\right)=
\frac{\partial V^{i}}{\partial x^{i}}+V^{i}\frac{\partial\ln \sigma}{\partial x^{i}}.
\label{Z4}$$
Then the \emph{Finsler-Laplacian} of $f$ can be defined by
$$\Delta f:=\textmd{div}(\nabla f)=\frac{1}{\sigma}\frac{\partial}{\partial x^{i}}\left(\sigma g^{ij}(\nabla f)f_j\right).\label{2.6}$$

\subsection{Anisotropic hypersurfaces}

Now let $\varphi :M^{n}\to (N^{n+1}, F)$ be an embedded hypersurface
of Finsler manifold $(N,F)$, where $M$ is an $n$-dimensional oriented smooth manifold. For simplicity, we will use $(x,y)\in
TM$, where $x=\varphi (u),~y=d\varphi v$ for $(u,v)\in TM$. Set
 $$\mathcal{V}(M)=\{\xi\in \varphi ^{-1}T^{*} {N}|~\xi (d\varphi (X))=0,~~\forall X\in \Gamma(TM)\},$$
which is called the {\it normal bundle} of $M$ in $(N,F)$.
 There exist exactly two unit normal vector fields
$\textbf{n}_{\pm}$ on $M$ such that $\textbf{n}_{\pm}={\mathcal L}^{-1}(\pm \nu )$, where $\nu\in\mathcal{V}(M)$ is a
unit 1-form. In general, $\textbf{n}_{+}$ is  not  necessarily $-\textbf{n}_{-}$ unless $F$ is reversible.

Let $\textbf{n}={\mathcal L}^{-1}(\nu )$ be a given normal vector of
$M$ in $(N,F)$, and put $ \hat{g}:=\varphi^{\ast}g_{\bf n}$. Then $(M, \hat{g})$ is a Riemannian manifold.  We also call $(M, \hat{g})$ an \emph{anisotropic hypersurface} of $(N,F)$ to distinguish  it from an isometric  immersion hypersurface of $(N,F)$. \begin{eqnarray*}
\hat{g}_{i j }={g}_{ \alpha  \beta }(\textbf{n})\varphi ^{ \alpha
}_i \varphi ^{\beta }_ j ,
\end{eqnarray*}
where
\begin{eqnarray*}
{x}^{ \alpha }=\varphi ^{ \alpha }(u),\quad  {y}^{ \alpha }=\varphi ^{ \alpha }_i
v^i ,\quad \varphi ^{ \alpha }_i =\frac {\partial \varphi ^{ \alpha
}}{\partial u^i }.
\end{eqnarray*}

Similar to the Riemannian case,  we have the \emph{Gauss-Codazzi equations} as follows (\cite{HYS})
\begin{equation}\label{2.14}\aligned{ D}^{\textbf{n}}_{X}Y&={ {\widehat{\nabla}}}_{X}Y+{h}(X,Y)\textbf{n},\\
D^{\textbf{n}}_{X}\textbf{n}&=-AX,~~~~~~~~~~~~~~~~~~~~~~\quad
X,Y\in \Gamma(TM) ,\endaligned\end{equation}
where \begin{equation}\label{2.15} h(X,Y):=\hat {g}(AX,Y)=\hat {g}(\textbf{n},D^{\textbf{n}}_{X}Y),
\quad X,Y\in
\Gamma(TM).\end{equation}
$\widehat{\nabla}$ is a torsion-free linear connection
on $M$ and satisfies
\begin{equation}\label{2.152}
(\widehat{\nabla}_X\hat g)(Y,Z)=2C_{\textbf{n}}(D^{\textbf{n}}_{X}\textbf{n},Y,Z)=-2C_{\textbf{n}}(AX,Y,Z),
\end{equation}
which shows that $\widehat{\nabla}$ is not the Levi-Civita connection on Riemannian manifold $(M,\hat{g})$.
$A:T_xM\rightarrow T_xM$ is
linear and self-adjoint with respect to $\hat g$, which is called the\emph{
Weingarten transformation} (or \emph{Shape operator}). So $h$ is bilinear and  $h(X,Y)=h(Y,X)$,
which is called the \emph{anisotropic second fundamental form} of $M$ in $(N,F)$.
 Moreover, we call the eigenvalues of $A$,
$k_1,k_2,\cdots,k_{n-1}$ , the \emph{anisotropic principal
curvatures} of $M$ with respect to $\textbf{n}$. If $A(X)=kX, \forall X\in \Gamma(TM)$, or equivalently, $k_1=k_2=\cdots=k_{n-1}$, we call $M$ an \textit{anisotropic totally umbilic} hypersurface of Finsler manifold $(N,F)$.

Define (\cite{Shen01, HYS})
\begin{equation}\label{2.16}\hat{H}:=\text {Tr} h=\sum_{a=1}^{n-1}k_{a}.\end{equation}  $\hat{H}$ is called
  the \emph{ anisotropic mean curvature} of $M$ in $(N,F)$.
 We remark that when $(N,F)$ is a Minkowski space, $\hat{H}$ is just the anisotropic mean curvature defined in many papers like \cite{Andrews} and \cite{Xia16}.

  Let
$d\mu=\sigma(x)dx^{1}\wedge\cdots\wedge dx^{n+1}$ be an arbitrary
volume form on $(N,F)$. The induced  volume form on $M$ determined
by $d\mu$ can be defined by
$$\aligned
d\mu_{\textbf{n}}={\sigma}({\varphi(u)})\varphi^*(i_{\bf n}(dx^1\wedge\cdots\wedge
dx^{n+1})),~~~ ~~~~~~u\in M,\endaligned$$
where $i_{\bf n}$
denotes the inner multiplication with respect to $\bf n$.
Consider a smooth variation $\varphi_t:M\to N$, $t\in (-\varepsilon
,\varepsilon )$, such that $\varphi_0=\varphi$.  From \cite{HYS}, we know that the first variation of the induced  volume along the direction $X$.
\begin{equation}\label{v1Proof7}\frac{d
{\text{Vol}}(t)}{dt}\Bigr|_{t=0}=-\int_{M}\mathcal{H}_{d\mu_{\textbf{n}}}(
{X})d\mu_{\textbf{n}},\end{equation}
 where
 $\mathcal{H}_{d\mu_{\textbf{n}}}$ is
called the \emph{$d\mu_{\textbf{n}}$-mean curvature form} of
$\varphi$. $H:=\mathcal{H}_{d\mu_{\textbf{n}}}(\textbf{n})$ is called the \emph{$d\mu_{\textbf{n}}$-mean curvature} of $M$ in $(N,F)$.
From \cite{Shen01}( also see \cite{HYS}), we know that \begin{equation}\label{2.170}H=\hat{H}+S(\textbf{n}),\end{equation} where $S(y)=y^i\frac{\delta}{\delta x^i}\left(
\ln\frac{\sqrt{\det(g_{ij})}}{\sigma(x)}\right)$, which called the S-curvature of $(N,F,d\mu)$. $M$ is called weakly convex (mean convex, resp.) if $h$ is nonnegative
definite ($\hat{H} \geq 0$ resp.). We say $h$ is strictly convex in $M$ if $h$ is
positive definite in $M$.

\noindent{\bf Lemma 2.2.1.} (\cite{HYS}) {\it In a Minkowski space $(N,F)$ with
the BH(resp. HT)-volume form $d\mu$, $\hat{H}=H$ and we have the fllowing Gauss-Codazzi equations:
\begin{eqnarray}\label{2.153}\aligned\hat{g}(\widehat{R}(X, Y)Z, W)=&{h}(Y, Z){h}(X, W)-{h}(X, Z){h}(Y, W),\\
(\widehat{\nabla}_{X}{h})(Y, Z)=&(\widehat{\nabla}_{Y}{h})(X, Z),
\endaligned\end{eqnarray}
for any $X, Y, Z, W \in \Gamma(TM)$, where $\hat{R}$ is the curvature tensor of the induced
connection $\widehat{\nabla}$ and $(\widehat{\nabla}_{X}{h})(Y, Z):=X\hat h(Y,Z)-\hat h({\hat{\nabla}}_XY,Z)-\hat h(Y,{\hat{\nabla}}_XZ)$.
}

\section{The Mean Curvature Flow in Finsler manifolds}

The geometric evolution equation \eqref{Fin1} is known as the mean curvature flow in Riemannian geometry. Similarly, we introduce the concept of the mean curvature
flow in the Finsler setting:

\noindent{\bf Definition 3.1.} {\it  Let $\varphi_{0} :(M^{n}, \hat{g})\to (N^{n+1}, F)$ be a smooth immersion of $M^{n}$ into a Finsler manifold $(N^{n+1},F)$. The mean curvature flow of $\varphi_{0}$ is a family of smooth immersions $\varphi_{t} :M\to (N^{n+1}, F)$
for $t\in[0, T)$ such that setting $\varphi(p,t)=\varphi_{t}(p)$ the map $\varphi :M\times[0, T)\to (N^{n+1}, F)$ is a smooth solution of the following system of PDE's
\begin{equation}\label{1Proof1}
\left\{\begin{array}{l}
\partial_{t}\varphi(p,t)={\hat{H}}(p,t)\textbf{n}(p,t),\\
\varphi(p,0)=\varphi_{0}(p),
\end{array}\right.
\end{equation}
where ${\hat{H}}(p,t)$ and $\textbf{n}(p,t)$ are respectively the anisotropic mean curvature and the inner pointing unit normal of the hypersurface $\varphi_{t}$ at the point $p\in M$.}

The following Gauss-Weingarten relations will be fundamental,
\begin{equation}\label{1Proof2}
\partial_{i}n^{\alpha}=-{h}^{j}_{i}\varphi^{\alpha}_j,\ \ \ \ \
\varphi^{\alpha}_{i|j}={h}_{ij}n^{\alpha}=\varphi^{\alpha}_{ij}+
\Gamma^{\alpha}_{\beta\gamma}\varphi^{\beta}_i\varphi^{\gamma}_j-\widehat{\Gamma}^{k}_{ij}\varphi^{\alpha}_k,
\end{equation}
where $\varphi^{\alpha}_{i|j}$ denotes the covariant derivative of $\varphi^{\alpha}_{i}$ and $\widehat{\Gamma}^{k}_{ij}$ is the Christoffel symbols with respects to $\hat{g}$.
Notice that, by these relations, it follows that
\begin{equation}\label{1Proof3}
\widehat{\Delta} \varphi=\hat{g}^{ij}\varphi_{i|j}
=\hat{g}^{ij}{h}_{ij}\textbf{n}={\hat{H}}\textbf{n}.
\end{equation}
Using equation \eqref{1Proof3}, this system can be rewritten in the appealing form
$$\partial_{t}\varphi=\widehat{\Delta}\varphi.$$

\noindent{\bf Proposition 3.1.}(Geometric Invariance under Tangential Perturbations). {\it
If a smooth family of immersions $\varphi_{t} :M\to (N^{n+1}, F)$ satisfies the system of PDE's
\begin{equation}\label{1Proof4}
\left\{\begin{array}{l}
\partial_{t}\varphi(p,t)={\hat{H}}(p,t)\textbf{n}(p,t)+X(p,t),\\
\varphi(p,0)=\varphi_{0}(p),
\end{array}\right.
\end{equation}
where $X$ is a time-dependent smooth vector field along $M$ such that $X(p, t)$ belongs
to $d\varphi_{t}(T_{p}M)$ for every $p\in M$ and every time $t\in [0, T)$, then, locally around any
point in space and time, there exists a family of reparametrizations (smoothly timedependent)
of the maps $\varphi_{t}$ which satisfies system \eqref{1Proof1}.\\
If the hypersurface M is compact, one can actually find uniquely a family of global
reparametrizations of the maps $\varphi_{t}$ as above for every $t\geq0$, leaving the initial
immersion $\varphi_{0}$ unmodified and satisfying system \eqref{1Proof1}.\\
Conversely, if a smooth family of moving hypersurfaces $\varphi_{t} :M\to (N^{n+1}, F)$
can be globally reparametrized for $t\geq0$ in order that it moves by mean curvature,
then the map $\varphi$ has to satisfy the system above for some time-dependent vector
field $X$ with $X(p, t) \in d\varphi_{t}(T_{p}M)$.}

\proof  First we assume that $M$ is compact; we will produce a smooth global
parametrization of the evolving sets in order to check the definition.

By the tangency hypothesis, the time-dependent vector field on $M$ given by
$Y(q, t)=-[d \varphi_{t}]^{-1}(X(q, t))$ is globally well defined and smooth.

Let $\Psi : M \times [0, T ) \to M$ be a smooth family of diffeomorphisms of $M$ with
$\Psi(p, 0) = p$ for every $p \in M$ and
\begin{equation}\label{1Proof5}
\partial_{t}\Psi(p, t)=Y(\Psi(p, t), t).
\end{equation}
This family exists, is unique and smooth, by the existence and uniqueness theorem
for ODE's on the compact manifold $M$.

Considering the reparametrizations
$\widetilde{\varphi}(p, t) = \varphi(\Psi(p, t), t)$, then
$$\aligned
\widetilde{\varphi}_i(p, t)=&\varphi_k(\Psi(p, t), t)\Psi^k_i,\\
\widetilde{\varphi}_{ij}(p, t)=&\varphi_{kl}(\Psi(p, t), t)\Psi^k_i\Psi^k_j+\varphi_{k}(\Psi(p, t), t)\Psi^k_{ij},\\
\widetilde{g}_{ij}(p, t)=&{g}_{kl}(\Psi(p, t), t)\Psi^k_i\Psi^l_j,\\
\widetilde{\textbf{n}}(p, t)=&\textbf{n}(\Psi(p, t), t),\\
\widetilde{h}_{ij}(p, t)=&{g}_{\widetilde{\textbf{n}}}(\widetilde{\textbf{n}},\widetilde{\varphi}_{ij}+\Gamma^{\alpha}_{\beta\gamma}\widetilde\varphi^{\beta}_i\widetilde\varphi^{\gamma}_j)\\
=&{g}_{{\textbf{n}}}({\textbf{n}},{\varphi}_{kl}+\Gamma^{\alpha}_{\beta\gamma}\varphi^{\beta}_k\varphi^{\gamma}_l)\Psi^k_i\Psi^l_j\\
=&{h}_{kl}(\Psi(p, t), t)\Psi^k_i\Psi^l_j,\\
\widetilde{H}(p, t)=&{\hat{H}}(\Psi(p, t), t).
\endaligned$$
From which, one has
$$\aligned
 \partial_{t}\widetilde{\varphi}(p, t)=&\partial_{t}\varphi(\Psi(p, t), t)+d \varphi(\Psi(p, t))(\partial_{t}\Psi(p, t))\\
=&{\hat{H}}(\Psi(p, t), t)\textbf{n}(\Psi(p, t), t)+X(\Psi(p, t), t)+d \varphi(\Psi(p, t))(Y(\Psi(p, t), t))\\
=&{\hat{H}}(\Psi(p, t), t)\textbf{n}(\Psi(p, t), t)+X(\Psi(p, t), t)\\
&-d\varphi(\Psi(p, t))([d \varphi(\Psi(p, t))]^{-1}X(\Psi(p, t), t))\\
=&{\hat{H}}(\Psi(p, t), t)\textbf{n}(\Psi(p, t), t)\\
=&\widetilde{H}(p, t)\widetilde{\textbf{n}}(p, t).
\endaligned$$
Hence, $\varphi$ satisfies system \eqref{1Proof1} and $\widetilde{\varphi}_{0} = \varphi_{0}$.

Conversely, this computation also shows that if
$\widetilde{\varphi}(p, t) = \varphi(\Psi(p, t), t)$ satisfies
system \eqref{1Proof1}, the family of diffeomorphisms $\Psi : M \times [0, T ) \to  M$ must solve
equation \eqref{1Proof5}, hence, it is unique if we assume $\Psi( \cdot, 0)= Id_{M}$ in order that the
map $\varphi_{0}$ is unmodified.

In the noncompact case, we have to work locally in space and time, solving the
above system of ODE's in some positive interval of time in an open subset $\Omega\subset M$
with compact closure, then obtaining a solution of system \eqref{1Proof1} in a possibly
smaller open subset of $\Omega$ and some interval of time.

Assume now that the reparametrized map $\widetilde{\varphi}(p, t) = \varphi(\Psi(p, t), t)$ is a mean
curvature flow. Differentiating, we get
$$\aligned
 \partial_{t}\widetilde{\varphi}(p, t)=&\partial_{t}\varphi(\Psi(p, t), t)+d \varphi(\Psi(p, t))(\partial_{t}\Psi(p, t))\\
=&\widetilde{H}(p, t)\widetilde{\textbf{n}}(p, t)\\
=&{\hat{H}}(\Psi(p, t), t)\textbf{n}(\Psi(p, t), t)
\endaligned$$
that is, $$\partial_{t}\varphi(q, t)={\hat{H}}(q, t)\textbf{n}(q, t)-d \varphi(q)(\partial_{t}\Psi(\Psi_{t}^{-1}(q), t)),$$
for every $q \in M$ and $t\in [0, T)$. Then, the last statement of the proposition follows
by setting $X(q, t)=-d \varphi_{t}(q)(\partial_{t}\Psi(\Psi_{t}^{-1}(q), t))$.
\endproof

\noindent{\bf Corollary 3.1.} {\it If a smooth family of hypersurfaces $\varphi_{t} = \varphi(\cdot , t)$ satisfies $g_{\textbf{n}}(\partial_{t}\varphi, \textbf{n})={\hat{H}}$, then it can be everywhere locally reparametrized to a mean curvature flow. If $M$ is compact, this can be done uniquely by global reparametrizations, without
modifying $\varphi_{0}$.}

We give a more geometric, alternative definition of the mean curvature flow as follow.

\noindent{\bf Definition 3.2.} {\it We still say that a family of smooth immersions $\varphi_{t} :M\to (N^{n+1}, F)$, for $t \in[0, T )$, a mean curvature flow if locally at every point, in space and time,
there exists a family of reparametrizations which satisfies system \eqref{1Proof1}.}

Let $\text{A}(t)$ is the area of closed hypersurface $\varphi_{t} :M\to (N^{n+1}, F)$ with respect to $d\mu_{\textbf{n}}$. Using \eqref{2.170} and \eqref{v1Proof7}, during the flow \eqref{1Proof1}, we have
\begin{equation}\label{v1Proof8}\frac{d
{\text{A}}(t)}{dt}=-\int_{M} (\hat{H}^{2}+\hat{H}S(\textbf{n})) d\mu_{\textbf{n}}.\end{equation}

\noindent{\bf Proposition 3.2.} {\it If a smooth family of closed hypersurfaces $\varphi_{t} = \varphi(\cdot , t)$ satisfies \eqref{1Proof1} and the Finsler manifold $(N^{n+1}, F,d\mu)$ has  vanishing ${\bf S}$-curvature, then the area of the hypersurfaces is nonincreasing i.e.,
$$\frac{d
{\text{A}}(t)}{dt}=-\int_{M} \hat{H}^{2} d\mu_{\textbf{n}}\leq0.$$
And thus we have the estimate
$$\int_0^{T_{max}}dt\int_{M} \hat{H}^{2} d\mu_{\textbf{n}}
\leq\text{A}(0).$$
}

\section{Examples}
\begin{expl} (Motion of Level Sets).\end{expl}
Level set representation $M_t=\{f(\varphi(p, t), t)=0\}$, where $f : N^{n+1} \times [0, T ) \to  \mathbb{R}$, $p \in M$.
Then we have
$$\textbf{n}=\frac{\nabla f}{F(\nabla f)},\ \ \ \ \ H=-\textmd{div}\Big(\frac{\nabla f}{F(\nabla f)}\Big)=\hat{H}+S(\textbf{n}).$$
We compute $$0=\frac{d (f(\varphi(p, t), t))}{dt}=g_{\nabla f}(\nabla f, \partial_{t}\varphi)+\partial_{t}f.$$
If $\varphi_{t}$ satisfies system \eqref{1Proof1}, then $g_{\nabla f}(\nabla f, \partial_{t}\varphi)=F(\nabla f)\hat{H}$ and thus $f$ satisfies
\begin{equation}\label{3.1.1}
\partial_{t}f=\Delta f-\frac{\nabla^2 f(\nabla f, \nabla f)}{F(\nabla f)^{2}}+S(\nabla f).
\end{equation}
Conversely, if we have a smooth function $f$ satisfying the above equation \eqref{3.1.1}, by Corollary 3.1,
every regular level set of $f( \cdot , t)$ is a hypersurface moving by mean curvature.

\begin{expl} (Isoparametric hypersurfaces).\end{expl}
Let  $f$ be an isoparametric function on a Finsler manifold $(N^{n+1}, F,d\mu)$ with constant {\bf S}-curvature $(n+1)cF$, that is, there is a smooth
function $ a (t)$ and a  continuous function $b (t)$
defined on $J=f(M_f)$ such that
\begin{equation}\label{3.1}\left\{\begin{aligned}
      &F(\nabla f)= a (f),\\
      &\Delta f= b (f).
\end{aligned}\right.\end{equation}
The isoparametric family of $f$ is $\{f=s\}$. Let $\varphi_{t} :M\to (N^{n+1}, F)$ be a family of smooth immersions satisfying
$$f(\varphi(p, t))=s(t).$$ Then from \cite{HYS}, the function  $\tilde{f}({x}, t)=f({x})-s(t)$ satisfies
$$\aligned
\partial_{t}\tilde{f}=&-s'(t),\\
\Delta \tilde{f}-\frac{\nabla^2\tilde{f}(\nabla \tilde{f}, \nabla \tilde{f})}{F(\nabla \tilde{f})^{2}}+S(\nabla \tilde{f})=&\Delta f-\frac{\nabla^2 f(\nabla f, \nabla f)}{F(\nabla f)^{2}}+c(n+1)F(\nabla \tilde{f})\\
=&b(s)-a(s)a'(s)+c(n+1)a(s)\\
=&a(s)\hat H.
\endaligned$$
If $\hat H\neq0$, we can set $t=\int \frac{1}{a(s)a'(s)-b(s)-c(n+1)a(s)}ds$. Then $\tilde{f}$ satisfies \eqref{3.1.1} and thus the isoparametric family $f=s(t)$ is a  mean curvature flow on Finsler manifolds. \cite{HYS}, \cite{HYS2017} and \cite{Xu} have given a great many examples of isoparametric hypersurfaces in some Finsler space forms like Minkowski
space with $BH$-volume (resp. $HT$-volume) form, Funk
space and Randers sphere with $BH$-volume form.
 \endproof

\begin{expl}  (Homothetically shrinking flows in a Minkowski space).\end{expl} A family
hypersurfaces that simply move by contraction during the evolution by
mean curvature is called a \emph{homothetically shrinking flow}.
 From\cite{HYS}, we know that a given hyperplane, a isoparametric family of  Minkowski hyperspheres(also called Wulff shapes) or $F^*$-Minkowski cylinders  are all special examples of homothetically shrinking flows in a Minkowski
space $(N^{n+1}, F)$.

\noindent{\bf Proposition 4.1.} {\it If an initial hypersurface $\varphi_{0} :M\to (N^{n+1}, F)$ satisfies ${\hat{H}}(p)+\lambda g_{\textbf{n}}(\varphi_{0}(p)-x_{0}, \textbf{n}_{0}(p))=0$ at every point $p \in M$ for some constant $\lambda > 0$ and $x_{0}\in N^{n+1}$, then it generates a homothetically shrinking mean curvature flow.\\
Conversely, if $\varphi :M\times[0, T)\to (N^{n+1}, F)$ is a homothetically shrinking mean curvature
flow around some point  $x_{0}\in N^{n+1}$ in a maximal time interval, then either ${\hat{H}}$ is identically zero or $${\hat{H}}(p, t)+\frac{g_{\textbf{n}}(\varphi(p, t)-x_{0}, \textbf{n}(p, t))}{2(T-t)}=0,$$
for every point $p \in M$ and time $t \in [0, T )$.}

\proof If the condition is satisfied, we consider the homothetically shrinking flow
$$\varphi(p, t)=x_{0}+\sqrt{1-2\lambda t}(\varphi_{0}(p)-x_{0}).$$
It is easy to prove that
$$\textbf{n}(p, t) = \textbf{n}_{0}(p),~~~~{\hat{H}}(p, t)=\frac{{\hat{H}}(p, 0)}{\sqrt{1-2\lambda t}}.
$$
Then we see that
$$g_{\textbf{n}}(\partial_{t}\varphi(p, t), \textbf{n}(p, t))=-\frac{\lambda g_{\textbf{n}}(\varphi_{0}(p)-x_{0}, \textbf{n}(p, t))}{\sqrt{1-2\lambda t}}=\frac{{\hat{H}}(p, 0)}{\sqrt{1-2\lambda t}}={\hat{H}}(p, t).$$
 Hence, by Corollary 3.1, this is a mean curvature flow of the
initial hypersurface $\varphi_{0}$, according to Definition 3.2.

Conversely, if the homothetically shrinking evolution
$$\varphi(p, t)=x_{0}+f(t)(\varphi_{0}(p)-x_{0})$$
is a mean curvature flow, for some positive smooth function $f :[0, T)\to \mathbb{R}$ with $f(0)=1$, $\lim_{t\to T}f(t)=0$ and  $f^{'}(t)\leq0$, by Corollary 3.1 we have $g_{\textbf{n}}(\partial_{t}\varphi, \textbf{n})={\hat{H}}$, hence
\begin{equation}\label{3.1}\aligned
{\hat{H}}(p, t)=&g_{\textbf{n}}(\partial_{t}\varphi(p, t), \textbf{n}(p, t))\\
=&f^{'}(t)g_{\textbf{n}}(\varphi_{0}(p)-x_{0}, \textbf{n}(p, t))\\
=&\frac{f^{'}(t)}{f(t)}g_{\textbf{n}}(\varphi(p, t)-x_{0}, \textbf{n}(p, t))\\
\endaligned\end{equation}
and
\begin{equation}\label{3.2}\aligned
{\hat{H}}(p, 0)=&f(t){\hat{H}}(p, t)\\
=&f(t)f^{'}(t)g_{\textbf{n}_0}(\varphi_0(p)-x_{0}, \textbf{n}_0(p)),
\endaligned\end{equation}
as $\textbf{n}(p, t)=\textbf{n}_{0}(p)$. If ${\hat{H}}\neq0$ at some point,  we have that $f(t)f^{'}(t)$ is equal to
some constant $C$ for every $t \in [0, T )$, combining (\ref{3.1}) and (\ref{3.2}). Hence, $f(t)=\sqrt{2Ct+1}$ as $f(0)=1$ and since $\lim_{t\to T}f(t)=0$, we conclude $f(t)=\sqrt{1-\frac{t}{T}}$. The thesis then follows from (\ref{3.1}). \endproof

\section{Short Time Existence of the Flow}
From this section, we suppose that $(N^{n+1}, F)$ is an $n+1$-dimensional Minkowski space.

\noindent{\bf Theorem 5.1.} {\it For any initial, smooth and compact hypersurface of $N^{n+1}$ given
by an immersion $\varphi_{0} :M\to (N^{n+1}, F)$, there exists a unique, smooth solution of system
\eqref{1Proof1} in some positive time interval. Moreover, the solution continuously depends on the initial immersion $\varphi_{0}$ in $C^{\infty}$.}

\proof  Let $\varphi_{0} :M\to (N^{n+1}, F)$ be a smooth immersion of a compact $n$-dimensional manifold.
For the moment we assume that this hypersurface is embedded, hence the inner pointing unit normal vector field $\overline{\textbf{n}}_{0}$ (about Euclidean metric)is globally defined and smooth.

We look for a smooth solution $\varphi :M\times[0, T)\to (N^{n+1}, F)$ of the parabolic problem \eqref{1Proof1} for some $T >0$.

Since we are interested in a solution for short time, we can forget about the immersion
condition ($d\varphi_{t}$ nonsingular) as it will follow automatically by the smoothness
of the solution and by the fact that $\varphi_{0}$ is a compact immersion, when $t$ is close to
zero.

Keeping in mind Proposition 3.1 and Corollary 3.1, if we find a smooth solution
of the problem
\begin{equation}\label{1Proof7}
\left\{\begin{array}{l}
g_{\textbf{n}}(\partial_{t}\varphi(p, t), \textbf{n}(p, t))={H}(p, t),\\
\varphi(p,0)=\varphi_{0}(p),
\end{array}\right.
\end{equation}
then we are done.

We consider the regular tubular neighborhood $\Omega=\{x\in N^{n+1}| d(x, \varphi_{0}(M))<\varepsilon\}$, which exists for $\varepsilon > 0$ small enough. By regular we mean that the map $\Psi :M\times(-\varepsilon, \varepsilon)\to \Omega$ defined as $\Psi(p, s)=\varphi_{0}(p)+s\overline{\textbf{n}}_{0}(p)$ is a diffeomorphism, where $\overline{\textbf{n}}_{0}$ is the unit normal vector of $\varphi_{0}$ with respect to the Euclidian metric in $N^{n+1}$.
Any small deformation of $\varphi_{0}(M)$ inside $\Omega$ can be represented as the graph of a
height function $f$ over $\varphi_{0}(M)$ and conversely, to any small function $f :M\to \mathbb{R}$ we can associate the hypersurface $M^{f}\subset \Omega$ given by $\varphi(p) = \varphi_{0}(p)+f(p)\overline{\textbf{n}}_{0}(p)$. We
want to compute now the equation for a smooth function $f$, time-dependent, in
order that $\varphi$ satisfies system \eqref{1Proof7}. Obviously, as $f(\cdot, 0)$ gives the hypersurface $\varphi_{0}$, we have $f(p, 0)= 0$ for every
$p\in M$.

First we compute the metric and the normal of the perturbed hypersurfaces.
It is easy to see that
\begin{equation}\label{1Proof8}\aligned
\varphi_{i}(p,0)=&\varphi_{0i}(p)+\overline{\textbf{n}}_{0}(p)f_{i}(p,t)-f(p,t) \overline{h}^{k}_{i}(p,0)\varphi_{0k}(p,0)\\
=&\varphi_{i}(p,f,\partial f),
\endaligned\end{equation}
where $\partial f=(f_1,\ldots,f_n)$.  Then the induced  metric and normal vector of $\varphi$ with respect to the Euclidian metric $\langle~ ,~\rangle$ in $N^{n+1}$ are
$$\aligned
\overline{g}_{ij}(p, t)=&\langle \varphi_{i}(p, t), \varphi_{j}(p, t)\rangle\\
=&\overline{g}_{ij}(p,0)-2f(p,t)\overline{h}_{ij}(p,0)+f^{2}(p,t)\overline{h}_{ik}\overline g^{kl}\overline{h}_{lj}(p,0)
+f_{i}(p,t)f_{j}(p,t)\\
=&\overline{g}_{ij}(p,f,\partial f),\\
\overline{\textbf{n}}=&\overline{\textbf{n}}_{0}-d\varphi \overline{\nabla}f,
\endaligned$$
where $\overline{\nabla}f$ is the gradient vector with respect to $\overline{g}$.
 Hence we know that
$$\aligned\nu=&\frac{\overline{\textbf{n}}}{F^{\ast}(\overline{\textbf{n}})}=\frac{\overline{\textbf{n}}_{0}-d\varphi \overline{\nabla}f}{F^{\ast}(\overline{\textbf{n}}_{0}-d\varphi \overline{\nabla}f)}=\nu(p,f,\partial f),\\
\textbf{n}=&{\mathcal L}^{-1}(\nu )=\textbf{n}(p,f,\partial f),\\
\hat{g}_{ij}(p, t)=&g_{\textbf{n}}(\varphi_{i}(p, t), \varphi_{i}(p, t))=\hat{g}_{ij}(p,f,\partial f).
\endaligned$$

Notice that the normal, the metric and thus its inverse depend only on first space
derivatives of the function $f$. Moreover, as $f(p, 0)=\overline{\nabla} f(p, 0) = 0$, everything is smooth and since $M$ is compact, when $t$ is small the denominator of the above
expression for the normal is uniformly bounded below away from zero (actually it is close to one).

From \eqref{1Proof8}, we have
$$\varphi_{ij}=\varphi_{0ij}+\overline{\textbf{n}}_{0}f_{ij}+f \overline{\textbf{n}}_{0ij}+f_{i}\overline{\textbf{n}}_{0j}
+f_{j}\overline{\textbf{n}}_{0i}.$$
Then, we find the second fundamental form,
$$\aligned
{h}_{ij}=&g_{\textbf{n}}(\textbf{n}, \varphi_{ij})=\nu(\varphi_{ij})
=\frac{1}{F^{\ast}(\overline{\textbf{n}})}\overline{h}_{ij}\\
=&\frac{\langle\overline{\textbf{n}}_{0}-d\varphi \overline{\nabla}f,\varphi_{ij}\rangle}{F^{\ast}(\overline{\textbf{n}}_{0}-d\varphi \overline{\nabla}f)} \\
=&\frac{1}{F^{\ast}(\overline{\textbf{n}}_{0}-d\varphi \overline{\nabla}f)}\Big(\overline{h}_{ij}(p,0)+f_{ij}+f\langle\overline{\textbf{\textbf{n}}}_{0}, \overline{\textbf{n}}_{0ij}\rangle-\langle d\varphi \overline{\nabla}f, \varphi_{0ij}\rangle\\
&-f_{i}\langle d\varphi \overline{\nabla}f, \overline{\textbf{n}}_{0j}\rangle-f_{j}\langle d\varphi \overline{\nabla}f, \overline{\textbf{n}}_{0i}\rangle-\langle d\varphi \overline{\nabla}f, f\overline{\textbf{n}}_{0ij}\rangle\Big)\\
=&\frac{1}{F^{\ast}(\overline{\textbf{n}}_{0}-d\varphi \overline{\nabla}f)}f_{ij}+P_{ij}(p,f,\partial f),
\endaligned$$
where $P_{ij}$ is a smooth form when $f$ and $\partial f$ are small, hence for $t$ small.

The mean curvature is then given by
$$\aligned
{H}=&\hat{g}^{ij}{h}_{ij}\\
=&\frac{1}{F^{\ast}(\overline{\textbf{n}}_{0}-d\varphi \overline{\nabla}f)}\hat{g}^{ij}f_{ij}+P(p,f,\partial f),
\endaligned$$
where $(\hat{g}^{ij})=(\hat{g}_{ij})^{-1}=(\hat{g}^{ij}(p,f,\partial f))$ and $P$ is a smooth function, assuming that $f$ and $\overline{\nabla} f$ are small.

We are finally ready to write down the condition $g_{\textbf{n}}(\partial_{t}\varphi(p, t), \textbf{n}(p, t))={H}(p, t)$ in terms of the function $f$,
$$\aligned
 (\partial_{t}f)g_{\textbf{n}}(\textbf{n}, \overline{\textbf{n}}_{0})&=g_{\textbf{n}}(\partial_{t}\varphi(p, t), \textbf{n}(p, t))\\
&={H}(p, t)\\
&=\frac{1}{F^{\ast}(\overline{\textbf{n}}_{0}-d\varphi \overline{\nabla}f)}\hat{g}^{ij}f_{ij}+P(p,f,\partial f).
\endaligned$$
Notice that $g_{\textbf{n}}(\textbf{n}, \overline{\textbf{n}}_{0})=\nu(\overline{\textbf{n}}_{0})=\frac{1}{F^{\ast}(\overline{\textbf{n}}_{0}-d\varphi \overline{\nabla}f)}$, we get
$$\aligned
 \partial_{t}f&=\hat{g}^{ij}f_{ij}+\frac{P(p, f, \partial f)}{\nu(\overline{\textbf{n}}_{0})}\\
 &=\hat{g}^{ij}f_{ij}+Q(p,f,\partial f),
\endaligned$$
where $Q(p, \cdot , \cdot )$ is a smooth function when its arguments are small.

Then, if the smooth function $f : M \times [0, T) \to \mathbb{R}$ solves the following partial
differential equation (before we had to deal with a system of PDE's)
\begin{equation}\label{1Proof9}
\left\{\begin{array}{l}
\partial_{t}f=\hat{g}^{ij}f_{ij}+Q(p,f,\partial f),\\
f(p,0)=0,
\end{array}\right.
\end{equation}
then $\varphi(p, t) = \varphi_{0}(p) + f(p, t)\overline{\textbf{n}}_{0}(p)$ is a solution of system \eqref{1Proof1} for the initial compact embedding $\varphi_{0}$, in a positive time interval.

Conversely, if we have a mean curvature flow $\varphi$ of $\varphi_{0}$, for small time the hypersurfaces
$\varphi_{t}$ are embedded in the tubular neighborhood $\Omega$ of $\varphi_{0}(M)$, then the function
$f(p, t)= \pi_{(-\varepsilon, \varepsilon)}[\Psi^{-1}(\varphi(p, t))]$ is smooth and $f(p, 0) = \pi_{(-\varepsilon, \varepsilon)}[\Psi^{-1}(\varphi(p, 0))]=\pi_{(-\varepsilon, \varepsilon)}[(p, 0)] = 0$, where $\pi_{(-\varepsilon, \varepsilon)}$ is the projection map on the second factor of
$M \times (-\varepsilon, \varepsilon)$, hence, by the above computations $f$ must solve problem \eqref{1Proof9}.

 PDE \eqref{1Proof9} is a quasilinear strictly parabolic equation. In particular it is not degenerate (in some sense, passing
to the height function $f$ we killed the degeneracy of systems \eqref{1Proof1} and \eqref{1Proof7})
hence, we can apply the (almost standard) theory of quasilinear parabolic PDE's.
The proof of a general theorem about existence, uniqueness and continuous dependence
of a solution for a class of problems including \eqref{1Proof9} can be found again
in \cite{Carlo11} (see Appendix A in \cite{Carlo11}).

Using the unique solution $f$ of problem \eqref{1Proof9} we consider the associated map
$\varphi(p) = \varphi_{0}(p)+f(p)\overline{\textbf{n}}_{0}(p)$, we possibly restrict the time interval in order that $\varphi_{t}$ are all immersions
and we apply Corollary 3.1 to reparametrize globally the hypersurfaces in a
unique way in order to get a solution of system \eqref{1Proof1}. This association is one-to-one,
as long as one remains inside the regular tubular neighborhood $\Omega$, hence, existence,
uniqueness, smoothness and dependence on the initial datum of a solution
of system \eqref{1Proof1} follows from the analogous result for quasilinear parabolic PDE's.

If the hypersurface is not embedded, that is, it has self-intersections, since
locally every immersion is an embedding, we only need a little bit more care in the
definition of the height function associated to a mean curvature flow (a regular
global tubular neighborhood is missing), in order to see that the correspondence
between a map $\varphi$ and its height function $f$ is still a bijection, then the same
argument gives the conclusion also in the nonembedded case.
\endproof

\section{Evolution of Geometric Quantities}

Set $\widehat C_{ijk}=C_{\alpha\beta\gamma}(\textbf{n})\varphi^{\alpha}_{i}\varphi^{\beta}_{j}\varphi^{\gamma}_{k},~\widehat C_{ijkl}=C_{\alpha\beta\gamma\delta}(\textbf{n})\varphi^{\alpha}_{i}\varphi^{\beta}_{j}\varphi^{\gamma}_{k}\varphi^{\delta}_{l}$,
where $C_{\alpha\beta\gamma\delta}(y)=\frac{\partial C_{\alpha\beta\gamma}(y)}{\partial y^{\delta}}$.

\noindent{\bf Lemma 6.1.} {\it
$$\aligned
&\hat{g}_{kl|i}=-2 \widehat C_{kls}h^{s}_{i},\\
&\hat{g}^{kl}_{|i}=2\widehat C^k_{pm}h^{m}_{i}\hat{g}^{pl},\\
& H_{|ij}=\widehat{\Delta} h_{ij}-Hh^{m}_{i}h_{jm}+|A|_{\hat{g}}^{2}h_{ij}
+2\widehat C^k_{pm}(h^{m}_{i}h^{p}_{k|j}+h^{m}_{j}h^{p}_{k|i})\\
&\qquad\qquad+2\widehat C^k_{lm}h^{m}_{i|j}h^l_{k}-2\widehat C^k_{mpq}h^{m}_{i}h^{p}_{j}h^q_{k}.
\endaligned$$}

\proof
Using \eqref{2.153},
we calculate the second covariant derivative of $H$ and obtain
$$\aligned
H_{|ij}=&(\hat{g}^{kl}h_{kl})_{|ij}\\
=&\hat{g}^{kl}h_{kl|ij}+\hat{g}^{kl}_{|ij}h_{kl}+\hat{g}^{kl}_{|i}h_{kl|j}+\hat{g}^{kl}_{|j}h_{kl|i}.
\endaligned$$
Hence,
$$\hat{g}^{kl}h_{kl|ij}=H_{|ij}-\hat{g}^{kl}_{|ij}h_{kl}-\hat{g}^{kl}_{|i}h_{kl|j}-\hat{g}^{kl}_{|j}h_{kl|i}.$$
We calculate the Laplacian of $h_{ij}$ and obtain
$$\aligned
\widehat{\Delta} h_{ij}=&\hat{g}^{kl}h_{ij|kl}=\hat{g}^{kl}h_{ik|jl}\\
=&\hat{g}^{kl}(h_{ik|lj}+\hat{R}_{ljkm}h^{m}_{i}+\hat{R}_{ljim}h^{m}_{k})\\
=&\hat{g}^{kl}(h_{kl|ij}+\hat{R}_{ljkm}h^{m}_{i}+\hat{R}_{ljim}h^{m}_{k})\\
=&\hat{g}^{kl}h_{kl|ij}+\hat{g}^{kl}h^{m}_{i}(h_{lk}h_{jm}-h_{lm}h_{jk})+\hat{g}^{kl}h^{m}_{k}(h_{mj}h_{li}-h_{ji}h_{lm})\\
=&\hat{g}^{kl}h_{kl|ij}+Hh^{m}_{i}h_{jm}-|A|_{\hat{g}}^{2}h_{ij}\\
=&Hh^{m}_{i}h_{jm}-|A|_{\hat{g}}^{2}h_{ij}+H_{|ij}-\hat{g}^{kl}_{|ij}h_{kl}-\hat{g}^{kl}_{|i}h_{kl|j}-\hat{g}^{kl}_{|j}h_{kl|i}.
\endaligned$$
Hence,
\begin{equation}\label{1Proof110}
H_{|ij}=\widehat{\Delta} h_{ij}-Hh^{m}_{i}h_{jm}+|A|_{\hat{g}}^{2}h_{ij}
+\hat{g}^{kl}_{|ij}h_{kl}+\hat{g}^{kl}_{|i}h_{kl|j}+\hat{g}^{kl}_{|j}h_{kl|i}.
\end{equation}
We calculate the first and second covariant derivative of $\hat{g}$ by using \eqref{2.152} and obtain
\begin{equation}\label{1Proof1101} \hat{g}_{kl|i}=-2C_{\alpha\beta\gamma}\varphi^{\alpha}_{k}\varphi^{\beta}_{l}\varphi^{\gamma}_{s}h^{s}_{i}
=-2\widehat C_{kls}h^{s}_{i}\end{equation}
and
\begin{equation}\label{1Proof1102} \aligned
\hat{g}_{ij|mn}=&(-2C_{\alpha\beta\gamma}\varphi^{\alpha}_{i}\varphi^{\beta}_{j}\varphi^{\gamma}_{q}h^{q}_{m})_{|n}\\
=&2C_{\alpha\beta\gamma\tau}\varphi^{\alpha}_{i}\varphi^{\beta}_{j}\varphi^{\gamma}_{q}
\varphi^{\tau}_{p}h^{q}_{m}h^{p}_{n}-2C_{\alpha\beta\gamma}\varphi^{\alpha}_{i}
\varphi^{\beta}_{j}\varphi^{\gamma}_{q}h^{q}_{m|n}\\
=&2\widehat C_{ijpq}h^{p}_{m}h^{q}_{n}-2\widehat C_{ijq}h^{q}_{m|n}.
\endaligned\end{equation}
It follows from \eqref{1Proof1101} and \eqref{1Proof1102} that
$$
\hat{g}^{kl}_{|i}=-\hat{g}^{ks}\hat{g}_{sp|i}\hat{g}^{pl}
=2C_{\alpha\beta\gamma}\varphi^{\alpha}_{s}\varphi^{\beta}_{p}\varphi^{\gamma}_{m}h^{m}_{i}\hat{g}^{ks}\hat{g}^{pl}
=2\widehat C_{spm}h^{m}_{i}\hat{g}^{ks}\hat{g}^{pl}.
$$
$$\aligned
\hat{g}^{kl}_{|ij}=&(2C_{\alpha\beta\gamma}\varphi^{\alpha}_{s}\varphi^{\beta}_{p}
\varphi^{\gamma}_{m}h^{m}_{i}\hat{g}^{ks}\hat{g}^{pl})_{|j}\\
=&2\hat{g}^{ks}_{|j}\hat{g}^{pl}C_{\alpha\beta\gamma}\varphi^{\alpha}_{s}\varphi^{\beta}_{p}\varphi^{\gamma}_{m}h^{m}_{i}
+2\hat{g}^{ks}\hat{g}^{pl}_{|j}C_{\alpha\beta\gamma}\varphi^{\alpha}_{s}\varphi^{\beta}_{p}\varphi^{\gamma}_{m}h^{m}_{i}\\
&-2C_{\alpha\beta\gamma\tau}\varphi^{\alpha}_{s}\varphi^{\beta}_{p}
\varphi^{\gamma}_{m}\varphi^{\tau}_{q}h^{q}_{j}h^{m}_{i}\hat{g}^{ks}\hat{g}^{pl}
+2C_{\alpha\beta\gamma}\varphi^{\alpha}_{s}\varphi^{\beta}_{p}\varphi^{\gamma}_{m}h^{m}_{i|j}\hat{g}^{ks}\hat{g}^{pl}\\
=&2(\hat{g}^{ks}_{|j}\hat{g}^{pl}+\hat{g}^{ks}\hat{g}^{pl}_{|j})\widehat C_{spm}h^{m}_{i}
-2\widehat C_{spmq}h^{q}_{j}h^{m}_{i}\hat{g}^{ks}\hat{g}^{pl}
+2\widehat C_{spm}h^{m}_{i|j}\hat{g}^{ks}\hat{g}^{pl}.
\endaligned$$
The fourth term of \eqref{1Proof110} can be written as
$$\aligned
\hat{g}^{kl}_{|ij}h_{kl}=&(2C_{\alpha\beta\gamma}\varphi^{\alpha}_{s}\varphi^{\beta}_{p}
\varphi^{\gamma}_{m}h^{m}_{i}\hat{g}^{ks}\hat{g}^{pl})_{|j}h_{kl}\\
=&2\hat{g}^{ks}_{|j}\hat{g}^{pl}C_{\alpha\beta\gamma}\varphi^{\alpha}_{s}\varphi^{\beta}_{p}\varphi^{\gamma}_{m}h^{m}_{i}h_{kl}
+2\hat{g}^{ks}\hat{g}^{pl}_{|j}C_{\alpha\beta\gamma}\varphi^{\alpha}_{s}\varphi^{\beta}_{p}\varphi^{\gamma}_{m}h^{m}_{i}h_{kl}\\
&-2C_{\alpha\beta\gamma\tau}\varphi^{\alpha}_{s}\varphi^{\beta}_{p}
\varphi^{\gamma}_{m}\varphi^{\tau}_{q}h^{q}_{j}h^{m}_{i}\hat{g}^{ks}\hat{g}^{pl}h_{kl}
+2C_{\alpha\beta\gamma}\varphi^{\alpha}_{s}\varphi^{\beta}_{p}\varphi^{\gamma}_{m}h^{m}_{i|j}\hat{g}^{ks}\hat{g}^{pl}h_{kl}.
\endaligned$$
The fifth term of \eqref{1Proof110} can be written as
$$\aligned
\hat{g}^{kl}_{|i}h_{kl|j}=&2C_{\alpha\beta\gamma}\varphi^{\alpha}_{s}
\varphi^{\beta}_{p}\varphi^{\gamma}_{m}h^{m}_{i}\hat{g}^{ks}\hat{g}^{pl}h_{kl|j}\\
=&2C_{\alpha\beta\gamma}\varphi^{\alpha}_{s}
\varphi^{\beta}_{p}\varphi^{\gamma}_{m}h^{m}_{i}\hat{g}^{ks}(h^{p}_{k|j}-\hat{g}^{pl}_{|j}h_{kl})\\
=&2C_{\alpha\beta\gamma}\varphi^{\alpha}_{s}
\varphi^{\beta}_{p}\varphi^{\gamma}_{m}h^{m}_{i}\hat{g}^{ks}h^{p}_{k|j}
-2C_{\alpha\beta\gamma}\varphi^{\alpha}_{s}
\varphi^{\beta}_{p}\varphi^{\gamma}_{m}h^{m}_{i}\hat{g}^{ks}\hat{g}^{pl}_{|j}h_{kl}
\endaligned$$
Similary, the sixth term of \eqref{1Proof110} can be written as
$$
\hat{g}^{kl}_{|j}h_{kl|i}
=2C_{\alpha\beta\gamma}\varphi^{\alpha}_{s}
\varphi^{\beta}_{p}\varphi^{\gamma}_{m}h^{m}_{j}\hat{g}^{ks}h^{p}_{k|i}
-2C_{\alpha\beta\gamma}\varphi^{\alpha}_{s}
\varphi^{\beta}_{p}\varphi^{\gamma}_{m}h^{m}_{j}\hat{g}^{ks}\hat{g}^{pl}_{|i}h_{kl}.
$$
Hence,
\begin{equation}\label{1Proof111}\aligned
 H_{|ij}=&\widehat{\Delta} h_{ij}-Hh^{m}_{i}h_{jm}+|A|_{\hat{g}}^{2}h_{ij}
+2\widehat C^k_{pm}(h^{m}_{i}h^{p}_{k|j}+h^{m}_{j}h^{p}_{k|i})\\
&+2\widehat C^k_{lm}h^{m}_{i|j}h^l_{k}-2\widehat C^k_{mpq}h^{m}_{i}h^{p}_{j}h^q_{k}.
\endaligned\end{equation}
\endproof

\noindent{\bf Proposition 6.2.}  {\it Under the mean curvature flow,
$$\aligned
&(1)\qquad \partial_{t}\textbf{n}=-d\varphi\widehat{\nabla}{H};\\
&(2)\qquad \partial_{t}\hat{g}_{ij}=-2{H}{h}_{ij}-2\widehat C(\partial_{i}, \partial_{j}, \widehat{\nabla}{H});\\
&(3)\qquad \partial_{t}\hat{g}^{ij}=2\hat{g}^{is}{h}_{s}^{j}{H}+2\hat{g}^{lj}\widehat C^i_{lk}H^k;\\
&(4)\qquad  \partial_{t}{H}=\widehat{\Delta} {H}+{H}|A|_{\hat{g}}^{2}+2tr_{\hat{g}}\widehat C(\cdot,A, \widehat{\nabla}{H});\\
&(5)\qquad\partial_{t}{h}_{ij}
=\widehat{\Delta} h_{ij}-2Hh^{m}_{i}h_{jm}+|A|_{\hat{g}}^{2}h_{ij}
+2\widehat C^k_{pm}(h^{m}_{i}h^{p}_{k|j}+h^{m}_{j}h^{p}_{k|i})\\
&\qquad\qquad\qquad+2\widehat C^k_{lm}h^{m}_{i|j}h^l_{k}-2\widehat C^k_{mpq}h^{m}_{i}h^{p}_{j}h^q_{k}.
\endaligned$$}

\proof  (1) Taking derivative of $\partial_{t}\varphi={H}\textbf{n}$ with respect to $x$ and
using  ${n}^{\alpha}_{i}=-h^{j}_{i}\varphi^{\alpha}_j$, we have
\begin{equation}\label{62Proof11}
(\partial_{t}\varphi)_{j}={H}_{j}\textbf{n}-{H}{h}^{k}_{j}\varphi_{k}.
\end{equation}
Taking derivative of ${n}^{\alpha}g_{\alpha\beta}(\textbf{n})\varphi^{\beta}_{i}=0$ with respect to $t$, we have
$$(\partial_{t}{n}^{\alpha})g_{\alpha\beta}(\textbf{n})\varphi^{\beta}_{i}
+2C_{\alpha\beta\gamma}(\partial_{t}{n}^{\gamma}){n}^{\alpha}\varphi^{\beta}_{i}
+{n}^{\alpha}g_{\alpha\beta}(\textbf{n})\partial_{t}\varphi^{\beta}_{i}=0.$$
From \eqref{62Proof11} and $C_{\alpha\beta\gamma}{n}^{\alpha}\varphi^{\beta}_{i}\partial_{t}{n}^{\gamma}=0$, we have
$$\aligned
(\partial_{t}{n}^{\alpha})g_{\alpha\beta}(\textbf{n})\varphi^{\beta}_{i}
=&-{n}^{\alpha}g_{\alpha\beta}(\textbf{n})\partial_{t}\varphi^{\beta}_{i}\\
=&-{n}^{\alpha}g_{\alpha\beta}(\textbf{n})({n}^{\beta}{H}_{i}
-{H}{h}^{k}_{i}\varphi^{\beta}_{k})\\
=&-{H}_{i}.
\endaligned$$
Using the above formulas, we have $\partial_{t}{n}^{\alpha}=-\hat{g}^{ij}\varphi^{\alpha}_{j}{H}_{i}.$
Hence, $$\partial_{t}\textbf{n}=-d\varphi\widehat{\nabla}{H}.$$

(2) Taking derivative of $\hat{g}_{i j }={g}_{ \alpha  \beta }(\textbf{n})\varphi ^{ \alpha
}_i \varphi ^{\beta }_ j$ with respect to $t$, we have
$$\aligned
\partial_{t}\hat{g}_{i j }
=&2(\partial_{t}\varphi ^{ \alpha}_i){g}_{ \alpha  \beta }(\textbf{n})\varphi ^{\beta}_j+2\varphi ^{ \alpha}_iC_{\alpha\beta\gamma}\varphi ^{\beta}_j(\partial_{t}{n}^{\gamma})\\
=&2{g}_{ \alpha  \beta }(\textbf{n})\varphi ^{\beta}_j({n}^{\alpha}{H}_{i}
-{H}{h}^{k}_{i}\varphi^{\alpha}_{k})-2\varphi ^{ \alpha}_iC_{\alpha\beta\gamma}\varphi ^{\beta}_j\hat{g}^{kl}\varphi^{\gamma}_{l}{H}_{k}\\
=&-2{H}{h}_{ij}-2\widehat{C}_{ijl}\hat{g}^{kl}{H}_{k}.
\endaligned$$

(3) Differentiating the formula $\hat{g}_{is}\hat{g}^{sj}= \delta^{j}_{i}$ we get
$$\aligned
\partial_{t}\hat{g}^{ij}=&-\hat{g}^{is}\partial_{t}\hat{g}_{sl}\hat{g}^{lj}\\
=&-\hat{g}^{is}\hat{g}^{lj}[-2{H}{h}_{sl}-2\widehat C(\partial_{s}, \partial_{l},\widehat{\nabla}{H})]\\
=&2\hat{g}^{is}{h}_{s}^{j}{H}+2\hat{g}^{is}\hat{g}^{lj} \widehat C(\partial_{s}, \partial_{l}, \widehat{\nabla}{H}).
\endaligned$$

(4) We calculate the evolution of $H$ and obtain
$$\aligned
\partial_{t}{H}=&\partial_{t}(\hat{g}^{ij}{h}_{ij})
=(\partial_{t}\hat{g}^{ij}){h}_{ij}+\hat{g}^{ij}\partial_{t}{h}_{ij}\\
=&{h}_{ij}[2\hat{g}^{is}{h}_{s}^{j}{H}
+2\hat{g}^{is}\hat{g}^{lj}\widehat C(\partial_{s}, \partial_{l}, \widehat{\nabla}{H})]
+\hat{g}^{ij}[\widehat{\nabla}^{2}_{ij}{H}-{H}{h}_{ik}\hat{g}^{kl}{h}_{lj}]\\
=&{H}|A|_{\hat{g}}^{2}+\widehat{\Delta} {H}+2tr_{\hat{g}}\widehat C_{n}(\cdot, A, \widehat{\nabla}{H}).
\endaligned$$

(5) From $\varphi^{\alpha}_{ij}=\widehat{\Gamma}^{k}_{ij}\varphi^{\alpha}_{k}+{h}_{ij}{n}^{\alpha}$, we have ${h}_{ij}=g_{\textbf{n}}(\textbf{n}, \varphi_{ij}).$
Using \eqref{62Proof11}, we have
$$
 (\partial_{t}\varphi)_{ij}=\textbf{n}{H}_{ij}-
 {H}_{i}{h}^{k}_{j}\varphi_{k}-{H}_{j}{h}^{k}_{i}\varphi_{k}
-{H}\partial_{j}{h}^{k}_{i}\varphi_{k}-{H}{h}^{k}_{i}\varphi_{kj}.
$$

Hence,
$$\aligned
\partial_{t}{h}_{ij}=&\partial_{t}g_{\textbf{n}}(\textbf{n}, \varphi_{ij})\\
=&g_{\textbf{n}}(\partial_{t}\textbf{n}, \varphi_{ij})+g_{\textbf{n}}(\textbf{n},  \partial_{t}\varphi_{ij})
+2C_{\alpha\beta\gamma}{n}^{\beta}(\partial_{t}{n}^{\gamma})\varphi_{ij}^{\alpha}\\
=&-g_{\textbf{n}}(d\varphi\widehat{\nabla} {H}, \varphi_{ij})+H_{ij}
-{H}{h}_{ki}{h}^{k}_{j}\\
=&-g_{\textbf{n}}({H}_{l}\varphi_{l}, \widehat{\Gamma}^{k}_{ij}\varphi_{k}+{h}_{ij}\textbf{n})+H_{ij}
-{H}{h}_{ki}{h}^{k}_{j}\\
=&-{H}_{l}\widehat{\Gamma}^{k}_{ij}\hat{g}_{kl}+H_{ij}
-{H}{h}_{ki}{h}^{k}_{j}\\
=& H_{|ij}-{H}{h}_{ik}\hat{g}^{kl}{h}_{lj}.
\endaligned$$
\endproof

\section{Consequences of Evolution Equations}

The main tool in order to obtain priori estimates is the maximum principle, in particular in the context of mean curvature flow. As in \cite{Carlo11, Andrews, Hamilton86}, we have the following maximum principles:

\noindent{\bf Lemma 7.1.}(\cite{Carlo11})  {\it Assume that $g(t)$, for $t \in [0, T)$, is a family of Riemannian metrics on a manifold $M$, with a possible boundary $\partial M$, such that the dependence on $t$ is smooth.

Let $u : M \times [0, T )\to \mathbb{R}$ be a smooth function satisfying
$$\partial_{t}u\geq \Delta_{g(t)}u+\langle X(p, u, \nabla u, t), \nabla u\rangle_{g(t)}+b(u)$$
where $X$ and $b$ are respectively a continuous vector field and a locally Lipschitz function in their arguments.

Then, suppose that for every $t \in [0, T)$ there exists a value $\delta > 0$ and a compact
subset $K \subset M \setminus\partial M$ such that at every time $t^{'} \in (t-\delta, t+\delta)\cap[0, T]$ the maximum
of $u( \cdot , t^{'})$ is attained at least at one point of $K$ (this is clearly true if $M$ is compact
without boundary).

Setting $u_{\min}(t)= \min_{p\in M} u(p, t)$ we have that the function $u_{min}$ is locally Lipschitz,
hence differentiable at almost every time $t \in [0, T )$ and at every differentiability time,
$$\frac{du_{\min}(t)}{dt}\geq b(u_{\min}(t)).$$
As a consequence, if $h : [0, T^{'} ) \to \mathbb{R}$ is a solution of the ODE
$$
\left\{\begin{array}{l}
h^{'}(t)=b(h(t)), \\
h(0)=u_{\min}(0),
\end{array}\right.
$$
for $T^{'} \leq T$, then $u \geq h$ in $M \times [0, T^{'})$.

Moreover, if $M$ is connected and at some time $\tau \in (0, T^{'})$ we have $u_{\min}(\tau) =
h(\tau)$, then $u = h$ in $M \times [0, \tau]$, that is, $u(\cdot , t)$ is constant in space.

}

\noindent{\bf Lemma 7.2.}(\cite{Andrews})  {\it Let $f\geq0$ be bounded and $C^{2}$ on $A\subset\mathbb{R}^{n}\times (0, T)$, such that
$$\frac{\partial f}{\partial t}\geq \alpha^{ij}f_{ij}+\beta^{i}f_{i}+\delta f,$$
where $\alpha^{ij}$ is non-negative definite and $C^{2}$, $\beta^{i}$ is $C^{1}$, and $\delta$ is bounded.
Suppose there exists $(x_{0}, t_{0})$ in the interior of $A$ such that $f(x_{0}, t_{0})=0$. Then $f(x, t)=0$
for every $(x, t)\in A$ for which there exists a piecewise smooth path $\gamma: [0, 1]\to A$ with $\gamma(0)=(x_{0}, t_{0})$, $\gamma(1)=(x, t)$, and
$$\gamma^{'}(s)=\sum\limits_{i,j=1}^{N}\omega_{i}(s)\alpha^{ij}\big|_{\gamma(s)}
\partial_{i}+r(s)\left(\left(\beta^{i}\big|_{\gamma(s)}-\partial_{j}\alpha^{ij}\big|_{\gamma(s)}
\right)(\partial_{i}-\partial_{t})\right)$$
with $r(s)\geq0$ and $\omega_{i}(s)\in \mathbb{R}$ for each $s\in[0, 1]$.
}

\noindent{\bf Lemma 7.3.}(\cite{Hamilton86})
Let $f(t)=\sup\{g(t,y):y\in Y\}$, where $Y$ is a
compact set. Then
$$\frac{d}{dt}f(t)\leq \sup \big(\frac{\partial}{\partial t}g(t, y):
y\in Y(t)\big), $$
where $Y(t)=\{y:g(t, y)=f(t)\}$.

Let us see some consequences of application of these maximum principles to evolution
equations for curvature. A hypersurface is mean convex if $H \geq 0$ everywhere. First we
will show that this property is preserved by the mean curvature flow.

\noindent{\bf Theorem 7.1.}(Preserving mean convexity) {\it Assume that the initial, compact hypersurface satisfies $H \geq 0$.
Then, under the mean curvature flow, the minimum of $H$ is increasing, hence $H$
is positive for every positive time. Moreover, $$H_{min}(t)\geq H_{min}(0)\left(1-\frac{2}{n}H_{min}^{2}(0)t\right)^{-\frac{1}{2}},$$
which gives an upper bound on the maximal existence time when $H_{min}(0)\neq0$:
$$T_{max}\leq\frac{n}{2}H_{min}^{-2}(0).$$
}

\proof  Arguing by contradiction, suppose that in an interval $(t_{0}, t_{1}) \subset \mathbb{R}^{+}$ we
have $H_{min}(t) < 0$ and $H_{min}(t_{0}) = 0$ ($H_{min}$ is obviously continuous in time and
$H_{min}(0) \geq 0$).

Let $|A|_{\hat{g}}^{2} \leq C$ in such an interval. Then
$$\partial_{t}{H}=\widehat{\Delta} {H}+2tr_{\hat{g}}C_{\textbf{n}}(d\varphi, d\varphi\circ A, d\varphi\widehat{\nabla}{H})+{H}|A|_{\hat{g}}^{2}$$
implies $$\partial_{t}{H}_{min}\geq C{H}_{min}$$
for almost every $t\in(t_{0}, t_{1})$.

Integrating this differential inequality in $[s, t]\subset (t_{0}, t_{1})$ we get
$H_{min}(t)\geq e^{C(t-s)}H_{min}(s)$, then sending $s \to t_{0}^{+}$ we conclude $H_{min}(t) \geq 0$ for every $t \in (t_{0}, t_{1})$ which is a contradiction.

Since then $H \geq 0$ we get
$$\aligned
\partial_{t}{H}=&\widehat{\Delta} {H}+2tr_{\hat{g}}C_{\textbf{n}}(d\varphi, d\varphi\circ A, d\varphi\widehat{\nabla}{H})+{H}|A|_{\hat{g}}^{2}\\
\geq&\widehat{\Delta} {H}+2tr_{\hat{g}}C_{\textbf{n}}(d\varphi, d\varphi\circ A, d\varphi\widehat{\nabla}{H})+\frac{{H}^{3}}{n}.
\endaligned$$
From Lemma 7.1, in $M_{t}$, $H(\cdot, t)\geq\varphi(t)$, where $\varphi(t)$ satisfying the ODE:
$$
\left\{\begin{array}{l}
\frac{\partial \varphi}{\partial t}=\frac{1}{n}\varphi^{3}, \\
\varphi(0)=H_{min}(0).
\end{array}\right.
$$
Clearly, the solution of the ODE is
$$\varphi(t)=H_{min}(0)\left(1-\frac{2}{n}H_{min}^{2}(0)t\right)^{-\frac{1}{2}}.$$
Then, if $H_{min}(0) = 0$ the ODE solution $\varphi(t)$ is always zero; so if at some positive
time $H_{min}(\tau) = 0$, we have, from Lemma 7.1, that $H( \cdot, \tau)$ is constant equal to zero on $M$, but there
are no compact hypersurfaces with zero mean curvature. Hence, $H_{min}$ is always
increasing during the flow and $H$ is positive on all $M$ at every positive time.

By the maximum principle, $$H_{min}(t)\geq\varphi(t)=H_{min}(0)\left(1-\frac{2}{n}H_{min}^{2}(0)t\right)^{-\frac{1}{2}}.$$
When $\big(1-\frac{2}{n}H_{min}^{2}(0)t\big)\rightarrow 0$, $H_{min}(t)\rightarrow \infty$, that is to say, $M_{t}$ has already blown up. However, $H$ is increasing and $H$ is positive for every positive time. Hence, from $1-\frac{2}{n}H_{min}^{2}(0)T\geq 0$, we have $T_{max}\leq\frac{n}{2}H_{min}^{-2}(0)$. \endproof

\noindent{\bf Teorem 7.2.}(Preserving convexity) {\it If $h_{ij}\geq0$ at $t=0$, then it remains so for $t>0$, that is, if the intial, the compact  hypersurface is convex, it remains convex under the mean curvature flow.}

\proof It is necessary to show that $h$ remains non-negative definite. A more general maximum principle is required here.

Let $\Theta$ be defined on $\{(p, v)\in TM: v\neq0\}$, the non-zero tangent bundle, by
\begin{equation}\label{eq:6a1}
\Theta(p, v, t)=\frac{h_{(p,t)}(v, v)}{\hat{g}_{(p,t)}(v, v)}
\end{equation}
for each $p\in M$ and $v\in T_{p}M\backslash \{0\}$. It suffices to show that $\Theta$ remains non-negative.
From Proposition 6.2, the evolution equation of $\Theta$ is
\begin{equation}\label{eq:6a2}\aligned
\partial_{t}\Theta=&\frac{1}{\hat{g}(v, v)}(\partial_{t}h_{_{ij}}-\Theta \partial_{t}\hat{g}_{ij})v^{i}v^{j}\\
=&\frac{1}{\hat{g}(v, v)}\{\widehat{\Delta} h_{ij}-2Hh^{m}_{i}h_{jm}+|A|_{\hat{g}}^{2}h_{ij}
+2\widehat C^k_{pm}(h^{m}_{i}h^{p}_{k|j}+h^{m}_{j}h^{p}_{k|i})\\
&+2\widehat C^k_{lm}h^{m}_{i|j}h^l_{k}-2\widehat C^k_{mpq}h^{m}_{i}h^{p}_{j}h^q_{k}
+2\Theta Hh_{ij}+2\Theta \widehat C(\partial_{i}, \partial_{j}, \widehat{\nabla}H)\}v^{i}v^{j}.
\endaligned\end{equation}
Given the coordinates $\{x^{i},  v^{i}\}$ on $TM$, one can observe that the pair $\{\frac{\delta}{\delta x^{i}},  \frac{\partial}{\partial v^{i}}\}$ forms a horizontal and vertical frame for $TTM$, where
$\frac{\delta}{\delta x^{i}}=\frac{\partial}{\partial x^{i}}-\widehat{\Gamma}^{k}_{ij}v^{j}\frac{\partial}{\partial v^{k}}$. Let $\{d x^{i},  \delta v^{i}\}$ denote the local frame dual to $\{\frac{\delta}{\delta x^{i}},  \frac{\partial}{\partial v^{i}}\}$, where $\delta v^{i}=d v^{i}+\widehat{\Gamma}^{k}_{ij}v^{j}dx^{j}$. Then we obtain a decomposition for $T(TM\backslash \{0\})$ and $T^{\ast}(TM\backslash \{0\})$,
$$T(TM\backslash \{0\})=\mathcal{H}TM\oplus \mathcal{V}TM,~~~~~~~T^{\ast}(TM\backslash \{0\})=\mathcal{H}^{\ast}TM\oplus \mathcal{V}^{\ast}TM,$$
where $$\mathcal{H}TM=span\{\frac{\delta}{\delta x^{i}}\},~~~~~~~\mathcal{V}TM=span\{\frac{\partial}{\partial v^{i}}\},$$
$$\mathcal{H}^{\ast}TM=span\{d x^{i}\},~~~~~~~\mathcal{V}^{\ast}TM=span\{\delta v^{i}\}.$$
Let $T=T^{j}_{i}\frac{\partial}{\partial x^{j}}\otimes dx^{i}$ be an arbitrary smooth local section of $\pi^{\ast}TM\otimes \pi^{\ast}T^{\ast}M$. The horizontal covariant derivative $T^{j}_{i|s}$ denotes
$$T^{j}_{i|s}=\frac{\delta T^{j}_{i}}{\delta x^{s}}+T^{k}_{i}\widehat{\Gamma}^{j}_{ks}-T^{j}_{k}\widehat{\Gamma}^{k}_{is}.$$
The vertical covariant derivative $T^{j}_{i;s}$ denotes
$$T^{j}_{i;s}=\frac{\partial T^{j}_{i}}{\partial v^{s}}.$$

From the above definite, we have
$$v^{l}_{|i}=\frac{\delta v^{l}}{\delta x^{i}}+\widehat{\Gamma}^{l}_{ij}v^{j}=0.$$
$$\aligned
h^{m}_{i|j}v^{i}v^{j}=&(\hat{g}^{mk}h_{ki})_{|j}v^{i}v^{j}\\
=&\hat{g}^{mk}_{|j}h_{k0}v^{j}+\hat{g}^{mk}h_{ij|k}v^{i}v^{j}\\
=&\hat{g}^{mk}_{|j}h_{k0}v^{j}+\hat{g}^{mk}(\Theta \hat{g}_{ij}v^{i}v^{j})_{|k}\\
=&\hat{g}^{mk}_{|j}h_{k0}v^{j}+\hat{g}^{mk}(\Theta_{|k} \hat{g}_{00}+\Theta \hat{g}_{00|k}),
\endaligned$$
where, for simplicity, we use the notations $h_{k0}=h_{kj}v^{j}$, $h_{00}=h_{ij}v^{i}v^{j}$, $\hat{g}_{k0}=\hat{g}_{kj}v^{j}$ and $\hat{g}_{00}=\hat{g}_{ij}v^{i}v^{j}$.

We also have
$$\Theta_{;k}=\frac{2}{\hat{g}_{00}}(h_{k0}-\Theta \hat{g}_{k0}).$$
$$\Theta_{|k}=\frac{1}{\hat{g}_{00}}(h_{00|k}-\Theta \hat{g}_{00|k}).$$
$$\aligned
\hat{g}^{ij}\Theta_{|i|j}
=&\hat{g}^{ij}\left(\frac{1}{\hat{g}_{00}}(h_{00|i}-\Theta \hat{g}_{00|i})\right)_{|j}\\
=&\frac{1}{\hat{g}_{00}}(\widehat{\Delta}h_{00}-2\hat{g}^{ij}\Theta_{|j}\hat{g}_{00|i}-\Theta \widehat{\Delta}\hat{g}_{00}).
\endaligned$$

$\Theta(p, v, t)=h_{(p,t)}(v, v)/\hat{g}_{(p,t)}(v, v)$ satisfies a degenerate
parabolic equation on the 2$n$-dimensional space $\{(p, v)\in TM: v\neq0\}$:
\begin{equation}\label{eq:61a2}
\partial_{t}\Theta=\hat{g}^{ij}\Theta_{|i|j}+a^{i}\Theta_{|i}+b^{i}\Theta_{;i}+d\Theta.\end{equation}
The coefficients are as follows:
$$\hat{g}^{ij}=\hat{g}^{ij}(x, t);$$
$$\aligned
a^{i}=a^{i}(x, v, t)=&\frac{2}{\hat{g}_{00}}\hat{g}^{ij}\hat{g}_{00|j}+2\widehat{C}^{k}_{lm}h^{l}_{k}\hat{g}^{mi}\\
=&\frac{-4}{\hat{g}_{00}}\hat{g}^{ij}\widehat{C}_{00s}h^{s}_{j}+2\widehat{C}^{k}_{lm}h^{i}_{k}\hat{g}^{mi};\endaligned$$
$$\aligned
b^{i}=b^{i}(x, v, t)=&-Hh^{i}_{j}v^{j}+2\widehat{C}^{k}_{pm}h^{p}_{k|j}v^{j}\hat{g}^{mi}
+\widehat{C}^{k}_{lm}h^{l}_{k}\hat{g}^{mi}_{|j}v^{j}-\widehat{C}^{k}_{mps}h^{p}_{j}h^{s}_{k}v^{j}\hat{g}^{mi}\\
=&-Hh^{i}_{j}v^{j}+2\widehat{C}^{k}_{pm}h^{p}_{k|j}v^{j}\hat{g}^{mi}
+2\widehat{C}^{k}_{lm}\widehat{C}_{spt}h^{l}_{k}h^{t}_{j}\hat{g}^{ms}\hat{g}^{ip}v^{j}-\widehat{C}^{k}_{mps}h^{p}_{j}h^{s}_{k}v^{j}\hat{g}^{mi};
\endaligned$$
$$\aligned
d=d(x,v, t)=&\frac{2}{\hat{g}_{00}}(2\widehat{C}^{k}_{pm}h^{p}_{k|j}v^{j}\hat{g}^{sm}\hat{g}_{s0}
+2\widehat{C}^{k}_{lm}\widehat{C}_{sqt}h^{l}_{k}h^{t}_{j}\hat{g}^{sm}\hat{g}^{pq}\hat{g}_{p0}v^{j}
-\widehat{C}^{k}_{mpq}h^{p}_{j}h^{q}_{k}v^{j}\hat{g}^{ms}\hat{g}_{s0}\\
&+\widehat{C}_{00pq}\hat{g}^{mn}h^{p}_{m}h^{q}_{n}-\widehat{C}_{00q}h^{q}_{m|n}\hat{g}^{mn}+\frac{\hat{g}_{00}}{2}|A|_{\hat{g}}^{2}
-2\widehat{C}^{k}_{lm}\widehat{C}_{00s}h^{l}_{k}h^{s}_{p}\hat{g}^{mp}
+\widehat{C}_{00l}\hat{g}^{kl}H_{k},
\endaligned$$
where, $\widehat{C}_{00s}=\widehat{C}_{ijs}v^{i}v^{j}$ and $\widehat{C}_{00pq}=\widehat{C}_{ijpq}v^{i}v^{j}$.
In particular, all of these coefficients are smooth and bounded for $t\in[0, T_{0}]\subset [0, T)$ on the unit sphere bundle.

Therefore at the minimum point, there exists a constant $C(T_{0})$ such that
$$\partial_{t}\Theta\geq-C(T_{0})\Theta,$$
so by Lemma 7.3 (also  see \cite{Hamilton86} , Lemma 3.5) the minimum of $\Theta$ decreases
no faster than exponentially, and convexity is never lost. \endproof

\noindent{\bf Theorem 7.3.}(Preserving strict convexity) {\it Let $\{M_{t}\}_{t>0}$ be a family of convex hypersurfaces moving by the flow \eqref{1Proof1}. Then $M_{t}$ is strictly convex for each $t > 0$.
}

\proof A version of the strong maximum principle will apply to show that the region
$$Z_{t}=\{p\in M: h_{(p, t)}(v, v)=0~~ for~~ some~~ v\in T_{p}M\backslash \{0\}\}$$
is open for each fixed $t > 0$. It follows that $Z_{t}$ is either empty or all of $M$. But $h_{ij}$ is positive definite wherever $M_{t}$ touches an enclosing ball, so the only possibility is that $Z_{t}$ is empty.

Let $p_{0}\in Z_{t}$ and $v_{0}\in T_{p_{0}}M$ such that $h_{(p_{0}, t)}(v_{0}, v_{0})=0$. Then there exists a neighborhood of $p_{0}$, $U(p_{0})$, such that $A\subset TM$, where $A\cong U(p_{0})\times \mathbb{R}^{n}$. Given the coordinates $\{x^{i},  y^{i}\}$ on $A$, we know that,  on the 2$n$-dimensional space $A$, \eqref{eq:61a2} can be rewritten as
$$\partial_{t}\Theta=\alpha^{IJ}\Theta_{IJ}+\tilde{a}^{I}\Theta_{I}+\tilde{d}\Theta,$$
where $\alpha^{IJ}=\left(
  \begin{array}{ccc}
\hat{g}^{ij} & \ast\\
    \ast & \ast\\
  \end{array}
\right) $, $\Theta_{IJ}=\big(\frac{\partial^{2}\Theta}{\partial x^{i}\partial x^{j}}, \frac{\partial^{2}\Theta}{\partial x^{i}\partial y^{j}}, \frac{\partial^{2}\Theta}{\partial y^{i}\partial y^{j}} \big)$, $\Theta_{I}=\big(\frac{\partial \Theta}{\partial x^{i}}, \frac{\partial \Theta}{\partial y^{i}}\big)$   and  all of these coefficients are smooth and bounded for $t\in[0, T_{0}]\subset [0, T)$ on the unit sphere bundle.

For every $p_{1}\in U(p_{0})$, there exists a curve $\sigma: [0, 1]\to U(p_{0})$ with $\sigma(s)=(x^{i}(s))$, where $\sigma(0)=p_{0}$ and $\sigma(1)=p_{1}$. Suppose that $\gamma(s):=(x^{i}(s), v^{i}(s))\in A$ such that $\gamma^{'}(s)=\dot x^{i}\frac{\partial}{\partial  x^{i}}+\dot v^{i}\frac{\partial}{\partial y^{i}}$, where $ v(s)$ is a vector field to be determined satisfying $v(0)=v_{0}$.
Since $v_{0}\neq0$, suppose $v_{0}^k\neq0$, we can choose $\omega_{I}(s)$ such that $\omega_{I}(s)\alpha^{I\overline{k}}(\sigma(s))=0$ and
$\omega_{I}(s)\alpha^{Ij}(\sigma(s))=\dot x^{j}(s)$ for any $s$, where $\overline{i}=n+i$. Let
\begin{equation}\label{Proof7}
\left\{\begin{array}{l}
\dot v^{j}=\omega_{I}\alpha^{I\overline{j}}\\
v(0)=v_{0},
\end{array}\right.
\end{equation}
Then $\gamma^{'}(s)=\omega_{I}\alpha^{IJ}\partial_{J}$ and $v^k(s)=v_{0}^k\neq0$, where $(\partial_{J})=(\frac{\partial}{\partial  x^{i}}, \frac{\partial}{\partial y^{i}})$.
Apply the strong maximum principle Lemma 7.2 with $f=\Theta$ and $r(s)$ identically zero on the domain $A$.   There exists a $v_{1}=v(1)\neq0$ such that $h_{(p_{1}, t)}(v_{1}, v_{1})=0$.  Therefore $p_{1} \in Z_{t}$, $U(p_{1}) \subset Z_{t}$, and $Z_{t}$ is open. \endproof

\section{Comparison Principle}

In order to study \eqref{1Proof1}, the following facts of convex hypersurfaces will be used.

Let $S^{n}_{F_{-}}(x_0;r):=\{x\in \mathbb{R}^{n}|F(-(x-x_0))=r\}$ be the inverse Minkowski hypersphere. For a convex hypersurface $M^{n}$ , we can also parametrize it as a graph over the unit inverse Minkowski hypersphere $S_{F_{-}}^{n}$. Let
\begin{equation}\label{com1}
\pi(x)=\frac{X(x)}{F(X(x))}: M^{n}\to S_{F_{-}}^{n},
\end{equation}
then we write the solution $M_{t}$ to \eqref{1Proof1} as a radial graph
\begin{equation}\label{com2}
X(x, t)=r(z,t)z:S_{F_{-}}^{n}\to \mathbb{R}^{n+1},
\end{equation}
where $r(z,t)=F(X(\pi^{-1}(z), t))$, $z\in S_{F_{-}}^{n}$ and $g_{\alpha\beta}(-z)z^{\alpha}z^{\beta}=1$.
The tangent vectors and the unit normal form of $M_t$ are given by
$$X^{\alpha}_{i}=r_{i}z^{\alpha}+rz^{\alpha}_{i},$$
$$\nu=\frac{\overline{\nu}}{F^{\ast}(\overline{\nu})},~~~\overline{\nu}_{\alpha}=g_{\alpha\beta}(-z)(-rz^{\beta}+r^{i}z^{\beta}_{i}),$$
where $F^{\ast}$ is the dual metric of $F$, $r^{i}=\overline{g}^{ij}r_{j}$ and $\overline{g}_{ij}=g_{\alpha\beta}(-z)z^{\alpha}_{i}z^{\beta}_{j}$. We have $$X^{\alpha}_{ij}=r_{ij}z^{\alpha}+rz^{\alpha}_{ij}+r_{i}z^{\alpha}_{j}+r_{j}z^{\alpha}_{i},$$
$$h_{ij}=\nu_{\alpha}X^{\alpha}_{ij}=\frac{1}{F^{\ast}(\overline{\nu})}
X^{\alpha}_{ij}g_{\alpha\beta}(-z)(-rz^{\beta}+r^{i}z^{\beta}_{i}).$$

 The mean curvature $H=\hat{g}^{ij}h_{ij}$ can be written in terms of $r(z,t)$ as follows:
\begin{equation}\label{com3}
H=-\frac{1}{F^{\ast}(\overline{\nu})}\hat{g}^{ij}(r\overline{\nabla}_{\partial_{i}}
\overline{\nabla}_{\partial_{j}}r-r^{2}\overline{g}_{ij}-2r_{i}r_{j}),
\end{equation}
where $(\hat{g}^{ij})=(\hat{g}_{ij})^{-1},~\hat{g}_{ij}=g_{\alpha\beta}(\textbf{n})X^{\alpha}_{i}X^{\beta}_{j}$, $\overline{\nabla}$ denotes the induced covariant derivative on $S_{F_{-}}^{n}$. Taking derivative of \eqref{com2} with respect to t, we have $r_{t}z = X_{t} = H\textbf{n}$. Since
$$H=\nu( X_{t} )=g_{-z}
\big( r_{t}z, \frac{-rz+r^{i}z_{i}}{F^{\ast}(\overline{\nu})}\big)=-\frac{r}{F^{\ast}(\overline{\nu})}r_{t},$$
thus
\begin{equation}\label{com4}
\partial_{t}r=\frac{-F^{\ast}(\overline{\nu})}{r}H.
\end{equation}
 By \eqref{com3} and \eqref{com4}, we get the evolution equation for $r(z, t)$,
\begin{equation}\label{com51}\partial_{t}r=\hat{g}^{ij}\big(  \overline{\nabla}_{\partial_{i}}
\overline{\nabla}_{\partial_{j}}r-\frac{2}{r}r_{i}r_{j}-r\overline{g}_{ij}\big).\end{equation}
If $H \geq 0$, then $r(z,t)=F(X(\pi^{-1}(z), t))$ is decreasing in
$[0, T)$ under the mean curvature flow equation \eqref{1Proof1}.

\noindent{\bf Proposition 8.1.} Assume that the initial, compact hypersurface satisfies $H \geq 0$, then $M_{t}$ is in fact contracting.

\noindent{\bf Proposition 8.2.} Let $S^{n}_{F_{-}}(O_{1};r_{1})$ and $S^{n}_{F_{-}}(O_{2};r_{2})$ be Minkowski hypersphere in $M^{n+1}$. If at $t=0$, $S^{n}_{F_{-}}(O_{1};r_{1})$ encloses the hypersurface $M$ and $S^{n}_{F_{-}}(O_{2};r_{2})$  is enclosed by the hypersurface $M$, then it remains so for $t>0$ and the distance between $S^{n}_{F_{-}}(O_{1};r_{1})$ and $M$(or $S^{n}_{F_{-}}(O_{2};r_{2})$  and $M$) is nondecreasing in time.

\proof
When $\overline{\nabla}r=0, \overline{\nu}=r \mathcal
L(-z)$. That is, $\textbf{n}=\frac{-z}{F(-z)}$, $g^{\alpha\beta}(\textbf{n})=g^{\alpha\beta}(-z)$,  $\hat{g}_{ij}=r^{2}\overline{g}_{ij}$ and $F^{\ast}(\overline{\nu})=r$. Thus, from \eqref{com51} we obtain that
$$\partial_{t}r\big|_{\overline{\nabla}r=0}=\overline{\Delta}r-\frac{n}{r}.$$

Note that the outer radius of $M$ is $r_{out}=r(z_1,t)$ and the inner radius of $M$ is $r_{in}=r(z_2,t)$. The distance
between $S^{n}_{F_{-}}(O_{1};r_{1})$ and $M$ (or $S^{n}_{F_{-}}(O_{2};r_{2})$ and $M$)is given by $r_{1}-r_{out}$ ( or $r_{in}-r_{2}$). we can conclude, as this analysis holds for all the pairs of points realizing the minimum, that
\begin{equation}\label{com5}
\partial_{t}(r_{1}-r_{out})=\overline{\Delta}(r_{1}-r)|_{z_1}-n\frac{r_{out}-r_{1}}{r_{1}r_{out}}\geq n\frac{r_{1}-r_{out}}{r_{1}r_{out}}.
\end{equation}
and
\begin{equation}\label{com6}
\partial_{t}(r_{in}-r_{2})=\overline{\Delta}(r-r_{2})|_{z_2}-n\frac{r_{2}-r_{in}}{r_{2}r_{in}}\geq n\frac{r_{in}-r_{2}}{r_{2}r_{in}}.
\end{equation}
Then $r_{1}-r_{out}$ ( or $r_{in}-r_{2}$) is nondecreasing in
$[0, T)$ under the mean curvature flow equation \eqref{1Proof1}. \endproof

FANQI ZENG\\
SCHOOL OF MATHEMATICAL SCIENCE
TONGJI UNIVERSITY\\
SHANGHAI 200092
P.R. CHINA\\
E-mail: fanzeng10$@$126.com\\

QUN HE\\
SCHOOL OF MATHEMATICAL SCIENCE
TONGJI UNIVERSITY\\
SHANGHAI 200092
P.R. CHINA\\
E-mail: hequn$@$tongji.edu.cn\\

BIN CHEN\\
SCHOOL OF MATHEMATICAL SCIENCE
TONGJI UNIVERSITY\\
SHANGHAI 200092
P.R. CHINA\\
E-mail: Chenbin$@$tongji.edu.cn\\

\end{document}